\documentclass[12pt]{article}
\usepackage{amsfonts}
\def\Bbb#1{\mathbb#1}
\newtheorem{theorem}{Theorem}

\newtheorem{proposition}[theorem]{Proposition}
\newtheorem{lemma}[theorem]{Lemma}
\newtheorem{corollary}[theorem]{Corollary}
\newtheorem{remark}[theorem]{Remark}

\newtheorem{example}[theorem]{Example}

\newcommand{\R}{\Bbb{R}}
\newcommand{\Q}{\Bbb{Q}}

\newcommand{\Sf}{\Bbb{S}}

\newcommand{\Ee}{\mathbb {E}}
\newcommand{\Les}{\Bbb{L}}

\newcommand{\Hy}{\Bbb{H}}

\newcommand{\spa}{\mbox{span}}

\newcommand{\rank}{\mbox{rank}}

\newcommand{\po}{{\hspace*{-1ex}}{\bf .  }}

\newcommand{\E}{{\cal E}}

\def\<{\langle}

\def\>{\rangle}
\def\a{\alpha}

\def\e{\epsilon}
\def\d{\partial}
\def\bea{\begin{eqnarray*} }
\def\eea{\end{eqnarray*} }
\def\be{\begin{equation} }
\def\ee{\end{equation} }

\def\nap{\nabla^\perp}

\def\proof{\noindent{\it Proof: }}
\def\qed{\ifhmode\unskip\nobreak\fi\ifmmode\ifinner\else
\hskip5 pt \fi\fi\hbox{\hskip5 pt \vrule width4 pt
height6 pt  depth1.5 pt \hskip 1pt }}

\begin{document}

\title{A decomposition theorem for\\ immersions of product manifolds.}
\author{Ruy Tojeiro\footnote{
Partially supported by CNPq grant 	311800/2009-2  and FAPESP grant
2011/21362-2.}}
\date{}
\maketitle

\begin{abstract}
We introduce polar metrics on a product manifold, which have product and warped product metrics as special cases. We prove a de Rham-type theorem  characterizing Riemannian manifolds that can be locally decomposed as a product manifold endowed with 
a polar metric. For  a product manifold endowed with a polar metric, our main result gives a complete description of all  its  isometric immersions into a space form whose second fundamental forms are adapetd to its product structure, in the sense that the tangent spaces to each factor are preserved by all shape operators. This is a far-reaching generalization of a basic decomposition theorem for isometric immersions of Riemannian products due to Moore as well as its extension by N\"olker to isometric immersions of warped products. 
\end{abstract}

\vspace{2ex}


\section{Introduction}

    A basic result due to Moore \cite{mo} states that an isometric immersion $f\colon\, M^n\to \R^N$ 
    of a Riemannian product $M^n=\Pi_{i=0}^r M_i$ of Riemannian manifolds is a product
    of isometric immersions $f_i\colon\, M_i\to \R^{m_i}$, $0\leq i\leq r$,
    with respect to some  orthogonal decomposition $\R^N=\oplus_{i=0}^r \R^{m_i}$, whenever the second fundamental form of 
    $f$ is adapted to the product structure of $M^n$, i.e.,  the tangent spaces to each factor are preserved 
    by all shape operators. 
    
    A significant generalization of Moore's theorem was obtained by N\"olker \cite{nol}. Under the same assumption on the 
    second fundamental form, he proved that an isometric immersion $f\colon\, M^n\to \Q_\epsilon^N$ into a simply connected space form of constant sectional curvature $\epsilon\in \{-1, 0, 1\}$
    of a warped product $M^n=M_0\times_{\rho}\Pi_{a=1}^r M_a$ of Riemannian manifolds splits as a warped  product
    of isometric immersions $f_i\colon\, M_i\to N_i$, $0\leq i\leq r$,
    with respect to some  warped product representation $\psi\colon\, N_0\times_{\sigma}\Pi_{a=1}^r N_a\to \Q_\epsilon^N$.

Our aim in this paper is to   prove a general decomposition theorem for immersions of product manifolds that is a far-reaching generalization of N\"olker's theorem.
Namely, we consider product manifolds $M^n=\Pi_{i=0}^r M_i$ endowed with a broad class of metrics that we call \emph{polar}.  These are  metrics   of the type
$$
g=\pi_0^*g_0+\sum_{a=1}^r \pi_a^*(g_a\circ \pi_{0})
$$
where $\pi_i\colon\,M\to M_i$, $0\leq i\leq r$, is the canonical projection, $g_0$ is a Riemannian metric on $M_0$ and $g_a$, $1\leq a\leq r$, are  smooth families of metrics  on $M_a$ parametrized by $M_0$.  

Polar metrics include product and warped product metrics as special cases. More generally, given Riemannian manifolds $M_0, \ldots, M_r$ with metrics $g_0, \ldots, g_r$, respectively, 
a metric $g$ on $M^n=\Pi_{i=0}^r M_i$ is called a \emph{quasi-warped product} of $g_0, \ldots, g_r$ with \emph{warping functions}
$\rho_a\colon M_0\times M_a\to \R_+$, $1\leq a\leq r$, if 
$$
g=\pi_0^*g_0+\sum_{a=1}^r (\rho_a\circ \pi_{0,a})^2\pi_a^*g_a,
$$
where $\pi_{0,a}\colon\,M^n\to M_0\times M_a$ is the projection. Warped product metrics 
$$
g=\pi_0^*g_0+\sum_{a=1}^r (\rho_a\circ \pi_{0})^2\pi_a^*g_a
$$
correspond to the case in which the warping functions $\rho_a$ are defined on $M_0$.  In particular, if all $\rho_a$ are identically one then $g$ is the Riemannian product of $g_0, \ldots, g_r$.

Therefore, a quasi-warped (respectively, warped, Riemannian) product  metric $g$ on  $M^n=\Pi_{i=0}^r M_i$ is a  polar metric 
 for which all metrics $g_a$ on $M_a$, $1\leq a\leq r$, are conformal  (respectively, homothetical, isometric) to a fixed Riemannian metric. 

We prove a de Rham type theorem characterizing Riemannian manifolds that can be locally decomposed as a product manifold endowed with 
a polar metric (see Theorem~\ref{cor:isodrh}). This can be regarded as a generalization of the local version of de Rham theorem, as well as its extensions by Hiepko \cite{hiepko} and the author \cite{to}, respectively, for warped product and quasi-warped product metrics. 

For  a product manifold endowed with a polar metric, our main result gives a complete description of all  its  isometric immersions into $\Q_\epsilon^N$ whose second fundamental forms are adapetd to its product structure (see Theorem \ref{thm:ptubes0} and Corollary \ref{cor:ptsphere} for the case of products with two factors, and Theorems~\ref{thm:ptubes} and \ref{thm:ptubeshyp} for the general case). The description relies on  the  concept of  partial tube  introduced by Carter and West~\cite{cw} (see Subsection~$3.2$). 
As a corollary, we obtain a decomposition theorem for isometric immersions into $\Q_\epsilon^N$ of quasi-warped product manifolds (see Theorems \ref{prop:n1c} and \ref{prop:n1e}). Then we show how N\"olker's theorem can be easily derived from the latter as a special case.

 As a consequence of our results, we obtain all representations of $\Q_\epsilon^N$ as a product manifold endowed with either a polar or a quasi-warped  product metric, i.e., all local isometries $\psi\colon\, \Pi_{i=0}^r N_i\to  \Q_\epsilon^N$ of such a product manifold into $\Q_\epsilon^N$. This generalizes N\"olker's classification of the so-called warped product representations of $\Q_\epsilon^N$, i.e., local isometries $\psi\colon\, \Pi_{i=0}^r N_i\to  \Q_\epsilon^N$ of a warped product manifold into $\Q_\epsilon^N$. 
 
 We conclude the paper with an application 
 to submanifolds with flat normal bundle of $\Q_\epsilon^N$. This includes a classification of surfaces with flat normal bundle of $\Q_\epsilon^N$ without umbilical points
 whose curvature lines of one family are geodesics.

    \section{Metrics on product manifolds}
    
    Our aim in this section is to characterize polar metrics on a product manifold $M^n=\Pi_{i=0}^r M_i$
in terms of the geometry of the product net of $M^n$. By a \emph {net} $\E=(E_i)_{i=0,\ldots,r}$ on a connected manifold $M^n$ 
we mean a splitting $TM=\oplus_{i=0}^r E_i$ of its tangent bundle by a family of  
integrable subbundles. The \emph{product net} on a  product manifold $M^n=\Pi_{i=0}^r M_i$ is the canonical net on $M^n$ determined by the tangent spaces to the factors. If $M^n$ is a Riemannian manifold and the subbundles of the net $\E$
are mutually orthogonal then $\E$ is said to be an \emph{orthogonal net}.

\subsection{Polar metrics}

    Given a product manifold
$M^n=\Pi_{i=0}^r M_i$, for
each $0\leq i\leq r$ we denote 
$$
M_{\perp_i}=M_0\times \cdots \times\hat{M}_i\times \cdots 
\times M_r,
$$ 
where  the hat over a term indicates that it is missing. The canonical projections of $M^n$ onto $M_i$ and $M_{\perp_i}$ are denoted by $\pi_i$ and $\pi_{\perp_i}$, respectively. 

For $\bar  p=(\bar p_0,\ldots, \bar p_r)\in M^n$, let  
 $\tau_i^{\bar p}\colon M_i\to M^n$  denotes the inclusion of $M_i$  into $M^n$ 
given  by 
$$
\tau_i^{\bar p}(p_i)=(\bar p_0, \ldots,p_i,\ldots, \bar p_r).
$$ 
Also, for any fixed $\bar p_i\in M_i$  the map $\mu_{\bar p_i}\colon M_{\perp_i}\to M^n$ stands for the inclusion of  $M_{\perp_i}$ into $M^n$ 
given  by 
$$
\mu_{\bar p_i}(p_0,\ldots,\hat p_i, \ldots, p_r)=(p_0, \ldots, \bar {p}_i, \ldots, p_r).
$$ 
Clearly, if $r=1$ then for any $\bar p=(\bar p_0, \bar p_1)\in M^n=M_0\times M_1$ we have $\tau_0^{\bar p}=\mu_{\bar p_1}$ and 
$\tau_1^{\bar p}=\mu_{\bar p_0}$.

The next lemma  characterizes the Riemannian metrics on a product manifold $M$
 for which the product net of $M$ is an orthogonal net.

\begin{lemma}\po\label{le:onet}  
Let $M^n=\Pi_{i=0}^r M_i$ be a product manifold, let $\E=(E_i)_{i=0,\ldots,r}$  
be its product net and let $g$ be a Riemannian metric on $M^n$. 
Then  $\E$ is an orthogonal net with respect to $g$ if and only if for
each $0\leq i\leq r$ 
there exists a smooth family  $g_i$ of metrics on $M_i$ parametrized by
$M_{\perp_i}$  such that 
\be\label{eq:metric0}
g=\sum_{i=0}^r \pi_i^*(g_i\circ \pi_{\perp_i}).
\ee
\end{lemma}
\proof  Given $X\in TM$, let  $X^i$ denote its $E_i$-component. Then, $\E$ is an orthogonal net with respect to $g$ if and only if 
\be\label{eq:metric0a} g(p)(X,Y)=\sum_{i=0}^rg(p)(X^i, Y^i)\ee
for any $p\in M^n$ and for all $X,Y\in T_pM$.

Given $p=(p_0,\ldots,  p_r)\in M^n$, set $p^i=\pi_{\perp_i}(p)$ and let  
 $g_i(p^i)=\tau_i^p{}^*g$ be the metric on
$M_i$ induced by $\tau_i^p$. 
Then, for any  $X,Y\in T_pM$ we have
\begin{eqnarray*}
g(p)(X^i, Y^i)&=&g(p)((\tau_i^p\circ
\pi_{i})_*{X^i}, (\tau_i^{p}\circ \pi_{i})_*{Y^i})=g_i( p^i)({\pi_i}_*{X^i},{\pi_i}_*{Y^i})\nonumber\vspace{1ex}\\&=&g_i( p^i)({\pi_i}_*{X},{\pi_i}_*{Y})=\pi_i^*(g_i(p^i))(X,Y),
\end{eqnarray*}
thus (\ref{eq:metric0a}) is equivalent to  (\ref{eq:metric0}).\vspace{1.5ex}\qed

To proceed  we need a well known characterization of the second fundamental form
 of an isometric immersion.
 
 Let $f\colon M^n\to \tilde{M}^m$ be an isometric immersion between Riemannian manifolds. A \emph{smooth
variation} of $f$ is a smooth mapping $F\colon I\times M^n\to\tilde{M}^m$, 
where $0\in I\subset\R$ is an open interval, such that  
$$
f_t=F(t,\cdot)\colon M^n\to \tilde{M}^m
$$ 
is an  immersion for any $t\in I$ and $f_0=f$. 

Let $\d/\d t$ denote the canonical vector field along the $I$ factor 
and set 
$$
T=F_*\d/\d t|_{t=0}
$$ 
regarded as a section of $f^*T\tilde{M}$. We say that $F$ is a \emph{normal
variation} if the variational vector field $T$ is everywhere normal to $f$. 
 
\begin{proposition}\po\label{prop:variation} 
Let $F\colon I\times M^n\to\tilde{M}^m$ be a smooth normal variation of $f$. 
Then 
$$
\frac{d}{d t}|_{t=0}\<{f_t}_*X, {f_t}_*Y\>=-2\<\a(X,Y),T\>\;\;\mbox{for
all}\;\;X,Y\in TM,
$$
where $\alpha\colon\,TM\times TM\to N_fM$ is the second fundamental form of $f$ 
with values in the normal bundle.
\end{proposition}
\proof  Consider the canonical extensions of  $X,Y$  over $I\times M^n$ and
note that 
$$
[X,\d/\d t]=0=[Y,\d/\d t].
$$ 
Then, 
\bea
\frac{d}{d t}\<{f_t}_*X, {f_t}_*Y\>
\!\!\!&=&\!\!\!\<\tilde{\nabla}_{\d/\d t}F_*X, F_*Y\>+\<F_*X,
\tilde{\nabla}_{\d/\d t}F_*Y\>\\
\!\!\!&=&\!\!\!\<\tilde{\nabla}_XF_*{\d/\d t}, F_*Y\>+\<F_*X,
\tilde{\nabla}_XF_*{\d/\d t}\>.
\eea
Hence, using the Weingarten formula we obtain
$$
\frac{d}{d t}|_{t=0}\<{f_t}_*X, {f_t}_*Y\>=-\<f_*A_TX,f_*Y\>-\<f_*X,f_*A_TY\>
=-2\<\a(X,Y),T\>.\,\,\,\qed
$$

We call a  metric $g$ on a product manifold $M_0\times M_1$ \emph{polar} if 
there exist 
a metric $g_0$ on $M_0$ and a smooth family $g_1$
of metrics  on $M_1$ parametrized by
$M_0$ such that 
$$
g=\pi_0^*g_0+\sum_{a=1}^r \pi_1^*(g_1\circ \pi_{0}).
$$
The terminology is explained by Example \ref{ex:polar}  in the next subsection. 

\begin{proposition}\po\label{prop:polarmetrica}  A metric $g$ on a product manifold $M=M_0\times M_1$ is polar if and
only if the product net $\E=(E_0, E_1)$ of $M$   is an orthogonal net  and $E_0$ is totally geodesic. 
\end{proposition}

Recall  that a  subbundle $E$ of the tangent bundle of  a Riemannian manifold  is said to be  \emph{totally geodesic} if 
$\nabla_XY\in \Gamma(E)$ for all $X, Y\in \Gamma(E)$. Here, for any vector subbundle $F$ of a vector bundle, we denote by $\Gamma(F)$ the space of its local smooth sections.\vspace{1ex}

\proof We already know from Lemma \ref{le:onet} that $\E$ is an orthogonal net with respect to $g$ if and only if 
$g$ is a metric of type (\ref{eq:metric0}). Given $p_1\in M_1$, let $g_0(p_1)$ be the metric on $M_0$ induced by the inclusion $\mu_{p_1}\colon\, M_0\to M_0\times M_1$. It follows from Proposition~\ref{prop:variation} that the metrics $g_0(p_1)$, $p_1\in M_1$,
all coincide if and only if $\mu_{p_1}$ is  totally geodesic  for any $p_1\in M_1$. Thus, a metric $g$ 
of type (\ref{eq:metric0}) is polar if and only if $E_0$ is totally geodesic.\qed\vspace{2ex}

As defined in the introduction, for a product manifold $M^n=\Pi_{i=0}^r M_i$ with arbitrarily many factors we call a metric $g$ on $M^n$ \emph{polar} if 
there exist a metric $g_0$ on $M_0$ and  smooth families $g_a$ 
of metrics on $M_a$ parametrized by $M_0$, $1\leq a\leq r$, such that 
$$
g=\pi_0^*g_0+\sum_{a=1}^r \pi_a^*(g_a\circ \pi_{0}).
$$

\begin{proposition}\po\label{le:productmetric}  Let $g$ be a Riemannian metric on a product manifold $M^n=\Pi_{i=0}^r M_i$.
Then the following assertions are equivalent:
\begin{itemize}
\item[$(i)$] $g$ is a polar metric;
\item[$(ii)$] The product net  $\E=(E_i)_{i=0,\ldots,r}$  of $M^n$ is an orthogonal net with respect to  $g$ such that  $E_a^\perp$ is totally geodesic  for all $1\leq a\leq r$;
\item[$(iii)$] $g$ is a polar metric on $M^n$, regarded as the product $M^n=M_{\perp_a}\times M_a$, for all $1\leq a\leq r$.
\end{itemize}
\end{proposition}
\proof The equivalence between $(ii)$ and $(iii)$ follows from Proposition \ref{prop:polarmetrica}. 
The condition in $(iii)$ is equivalent to $g$ being a metric of type (\ref{eq:metric0}) with the property that
the metrics induced on the leaves of the
product foliation of $M^n$ correspondent to $M_{\perp_a}$  coincide. Since a leaf
of the product foliation of $M^n$ correspondent to $M_{i}$ is contained in the
leaf correspondent to $M_{\perp_j}$ for any $j\neq i$, it follows that this property is satisfied
if and only if  the
metrics induced on the leaves of the product foliation of $M^n$ correspondent to
$M_{i}$, $i\neq 0$, do not depend on $j$ for every $1\leq j\leq r$ with $j\neq
i$, whereas all those correspondent to $M_0$ 
coincide. These are the conditions for $g$ to be a polar metric. Thus $(i)$ and $(iii)$ are equivalent.\vspace{1ex}\qed

For later use, we state a characterization  obtained in \cite{to} of the additional geometric properties that the product net of a product manifold must have with respect to a Riemannian metric $g$ 
  in order that $g$  be a quasi-warped  product metric.

 A subbundle $E$ of the tangent bundle of a Riemannian manifold  is said to be  \emph{umbilical} if there
exists a vector field $\eta$ in $E^\perp$, called the \emph{mean
curvature normal} of $E$, such that
$$\langle \nabla_XY,Z\rangle=\langle X,Y\rangle \langle\eta,
Z\rangle\,\,\,\,\mbox{for all}\,\,\,X, Y\in \Gamma(E), \,\,Z\in \Gamma(E^\perp).$$ 

\begin{proposition}\po\label{cor:prodwarp}  A Riemannian metric $g$ on a product manifold $M^n=\Pi_{i=0}^r M_i$ is a quasi-warped product metric if and only if the product net  $\E=(E_i)_{i=0,\ldots,r}$ of $M^n$ is an orthogonal net such that $E_a$ is umbilical and $E_a^\perp$ is totally geodesic 
for all $1\leq a\leq r$.
\end{proposition}

Warped product metrics (respectively, product metrics)  on $M^n=\Pi_{i=0}^r M_i$ are similarly characterized by 
replacing the word "umbilical" by "spherical" (respectively, "totally geodesic") in the previous statement (see \cite{hiepko}). 
That a subbundle $E$ of the tangent bundle of a Riemannian manifold  is \emph{spherical} means that it is umbilical and, in
addition, its mean
curvature normal $\eta$ satisfies
$$\langle \nabla_X\eta,Z\rangle=0,\,\,\,\mbox{for all}\,\,\,X\in \Gamma(E),\,\, Z\in \Gamma(E^\perp).$$

\subsection{A de Rham type theorem}

In order to state a precise de Rham-type characterization of Riemannian manifolds that can be locally decomposed as a product manifold endowed with a polar metric,  we introduce some further terminology (following \cite{rs}).

A
$C^\infty$-map $\psi\colon\;M\to N$ between two {\it netted
manifolds\/}  $(M,\mathcal{E})$, $(N,{\mathcal{F}})$, that is,
$C^\infty$-manifolds $M,N$ equipped with nets $\mathcal{E}=(E_i)_{i=0,\ldots,k}$ and $\mathcal{F}=(F _i)_{i=0,\ldots,k}$,
 respectively, is called a {\it net morphism\/} if $\psi_*E_i(p)\subset F_i(\psi(p))$
 for all $p\in M$, $0\leq i\leq k$,
or equivalently, if for any $p\in M$ the restriction
$\psi|_{L_i^{\mathcal{E}}(p)}$ to the leaf of $E_i$ through $p$ is a
$C^\infty$-map into the leaf $L_i^{\mathcal{F}}(\psi(p))$ of $F_i$
through $\psi(p)$. The net morphism $\psi$ is said to be a {\it
net isomorphism\/} if, in addition, it is a diffeomorphism and
$\psi^{-1}$ is also a net morphism. A net $\mathcal{E}$ on $M$ is
said to be {\em locally decomposable\/} if for every point $p\in
M$ there exist a neighborhood $U$ of $p$ in $M$ and a net
isomorphism $\psi$ from $(U,\mathcal{E}|_U)$ onto a product manifold
$\Pi_{i=0}^k M_i$. The map $\psi^{-1}\colon\;\Pi_{i=0}^k M_i\to U$
is called a {\em product representation\/} of $(U,\mathcal{E}|_U)$.

\begin{example}\po\label{ex:polar}{\em Let $G$ be a closed subgroup of the isometry group of a complete Riemannian manifold $M$ acting polarly on $M$. This means that there exist complete submanifolds of $M$, called \emph{sections}, 
that meet every orbit orthogonally. This turns out to be equivalent to the distribution $\mathcal H$ of normal spaces to 
the maximal dimensional orbits to be  integrable, in which case
$\mathcal H$ is easily seen to be totally geodesic. 

Fix a $G$-regular point $\bar p\in M$ and let $\Sigma$ be the section through $\bar p$. 
Then the map $\psi\colon\, \bar M=\Sigma\times G\bar p\to M_r$, onto the regular part of $M$,  given by $$\psi(p_0, g\bar p)=g(p_0),$$ 
 is a local diffeomorphism which is a net morphism with respect to 
the product net of $\Sigma\times G\bar p$ and the orthogonal net $({\cal H}, {\cal V})$ on $M_r$ given by its horizontal and vertical distributions ${\cal H}$ and ${\cal V}={\cal H}^\perp$, respectively. Since $\cal H$ is totally geodesic, it follows from Proposition \ref{prop:polarmetrica} that the metric on $\bar M$ induced by $\psi$ is polar.}
\end{example}

\begin{theorem}\po \label{cor:isodrh} 
Let $M$  be a Riemannian manifold with an orthogonal net $\mathcal{E}=(E_i)_{i=0,\ldots,k}$ such that 
$E_a^\perp$ is  totally geodesic   for each $1\leq a\leq r$. Then, for every point $p\in M$ there exists a
local product representation $\psi\colon\;\Pi_{i=0}^k M_i\to U$ of $\mathcal{E}$
with $p\in U\subset M$, which is an isometry with respect to a 
 polar metric  on $\Pi_{i=0}^k M_i$.
\end{theorem}

 Theorem \ref{cor:isodrh} is a consequence of  Proposition \ref{le:productmetric}
    and the following basic criterion for local decomposability of a net on an
    arbitrary $C^\infty$-manifold (cf. Theorem~$1$ of \cite{rs}).
    
\begin{proposition}\po \label{prop:dec} $\mbox{{\bf \cite{rs}}}$
A net $\mathcal{E}=(E_i)_{i=0,\ldots,k}$ on a $C^\infty$-manifold is
locally decomposable if and only if $E_i^{\perp}:=\oplus_{j\neq i}
E_i$ is integrable for $i=0,\ldots, k$.
\end{proposition}

\section{Immersions of product manifolds}

In this section, we prove a general decomposition theorem for immersions into space forms of product manifolds endowed 
with polar  metrics  (see Theorem \ref{thm:ptubes0} and Corollary~\ref{cor:ptsphere}  below for the case of products with two factors and 
Theorems \ref{thm:ptubes} and \ref{thm:ptubeshyp} for the general case).

For the convenience of the reader, we have included two preliminary subsections. In the first one we introduce
basic definitions on products of  immersions and state Moore's decomposition theorem \cite{mo}, as well as its extension due to Molzan \cite{mol} for isometric immersions into the sphere and hyperbolic space. 
The second one is
devoted to a special case of the notion of a \emph{partial tube} introduced by Carter and West \cite{cw}, on 
which our result relies.

\subsection[Products of  immersions]{Products of immersions}

Given  immersions $f_i\colon M_i^{n_i}\to\R^{m_i}$, $0\leq i\leq r$, 
their \emph{product} is the map
$$
f=f_0\times\cdots\times f_r\colon M^n=\Pi_{i=0}^r M_i^{n_i}\to\R^N
=\Pi_{i=0}^r\R^{m_i}
$$ 
given by $$f(p_0,\ldots,  p_r)=(f_0(p_0),\ldots,  f_r(p_r)).$$

If $f_i(M_i^{n_i})$ is contained in a hypersphere $\Sf^{m_i-1}(R_i)$ of radius $R_i$ centered at the origin, $0\leq i\leq r$, then 
$f(M^n)$ is contained in the hypersphere $\Sf^{N-1}(R)$ centered at the origin of radius $R$ given by $R^2=\sum_{i=0}^r R_i^2$.
If $f$ is regarded as an  immersion into $\Sf^{N-1}(R)$ and $f_i$ as an immersion into $\Sf^{m_i-1}(R_i)$
for $0\leq i\leq r$, then $f$ is called the \emph{extrinsic product} of $f_0, \ldots, f_r$.

In a similar way, let $\Hy^{N-1}(R)$ denote the hyperbolic space
$$\Hy^{N-1}(R)=\{X=(x_1, \ldots, x_n)\in \Les^N\,:\,\<X, X\>=-\frac{1}{R^2},\,\,x_1>0\}$$
of constant sectional curvature $-1/R^2$, where $\Les^N$ is Lorentz space of dimension $N$. Given an orthogonal decomposition $$\Les^N=\Les^{m_0}\times \Pi_{a=1}^r\R^{m_a}$$ one can define the extrinsic product $$f\colon\, M^n=\Pi_{i=0}^r M_i^{n_i}\to \Hy^{N-1}(R) \subset \Les^N$$   of  immersions $f_0\colon\, M_0^{n_0}\to \Hy^{m_0-1}(R_0) \subset \Les^{m_0}$ and $f_a\colon\, M_a^{n_a}\to \Sf^{m_a-1}\subset \R^{m_a}$, $1\leq a\leq r$, with 
$-R^2=-R_0^2+\sum_{a=1}^r R_a^2$.

A few elementary properties of a product 
$$
f=f_0\times\cdots\times f_r\colon M^n=\Pi_{i=0}^r M_i^{n_i}\to\R^N
=\Pi_{i=0}^r\R^{m_i}
$$ of  immersions $f_i\colon M_i^{n_i}\to\R^{m_i}$, $0\leq i\leq r$, are collected in the next proposition.

\begin{proposition}\po \label{prop:prod} The
following holds:
\begin{itemize}
\item[(i)] The differential of $f$ at any $p=(p_0,\ldots,  p_r)\in M^n$ is given by
$$
f_*\tau_i^p\!{}_*\,X_i=f_i{}_*\,X_i\,\,\,\mbox{for any}\,\,X_i\in T_{p_i}M_i,\,\,\,0\leq i\leq r.
$$
\item[(ii)] The map $f$ is an immersion whose  induced metric 
is the Riemannian product of the Riemannian metrics on $M_i$ induced by $f_i$, $0\leq i\leq r$. 
\item[(iii)] The normal space of  $f$ at $p$ is
$$
N_fM(p)=\oplus_{i=0}^r N_{f_i}M_i(p_i).
$$
\item[(iv)] The second fundamental form of $f$ at $p$ is given by
$$\alpha_f(\tau_i^p\!{}_*\,X_i,\tau_j^{p}\!{}_*\,X_j)=0\,\,\,\mbox{for any}\,\,\,X_i\in T_{p_i}M_i,\, X_j\in T_{p_j}M_j, \,\,i\neq j,$$
$$\alpha_f(\tau_i^p\!{}_*\,X_i,\tau_i^p\!{}_*\,Y_i)=\alpha_{f_i}(X_i, Y_i)\,\,\,\mbox{for all}\,\,\,X_i, Y_i\in T_{p_i}M_i,$$
or equivalently,
\be\label{eq:so} A^f_\xi\tau_i^p\!{}_*=\tau_i^p\!{}_*A^{f_i}_{\xi_i},\,\,\,0\leq i\leq r, \,\,\,\mbox{for all}\,\,\,\xi=\xi_0+\ldots +\xi_r\in N_fM.\ee
\end{itemize}
\end{proposition}

 The second fundamental form $\alpha$ of an isometric immersion 
$f\colon M^n\to \R^{N}$  is said to be \emph{adapted to the
net} $\E=(E_i)_{i=0,\ldots,r}$ on $M^n$ if $\alpha(E_i, E_j)=0$ for $0\leq i\neq j\leq r$.
It follows from part $(iv)$ of  Proposition \ref{prop:prod}  that the second fundamental form of a product
$$f=f_0\times\cdots\times f_r\colon M^n=\Pi_{i=0}^r M_i^{n_i}\to\R^N
=\Pi_{i=0}^r\R^{m_i}$$  of isometric immersions $f_i\colon M_i^{n_i}\to\R^{m_i}$, $0\leq i\leq r$, is 
adapted to the product net of $M^n$. The next  useful result due to Moore \cite{mo} shows that products of isometric immersions are  characterized by this property among isometric immersions of Riemannian products.

\begin{theorem} \label{thm:moore2} Let $f\colon M^n=\Pi_{i=0}^r M_i^{n_i}\to \R^{N}$ be
an isometric immersion of a Riemannian product  whose second
fundamental form is adapted to the product net of $M^n$.  Then there exist  an orthogonal decomposition
$\R^N=\Pi_{i=0}^r \R^{m_i}$ and isometric immersions $f_i\colon M_i\to \R^{m_i}$, $0\leq i\leq r$, such 
that $f$ is the product of $f_0,\ldots,f_r$.
\end{theorem} 

Moore's theorem has been generalized by Molzan \cite{mol} for isometric immersions into the  sphere and hyperbolic space.
For instance, in the latter case his result can be stated as follows. 

\begin{theorem} \label{thm:molzan} Let $f\colon M^n=\Pi_{i=0}^r M_i^{n_i}\to \Hy^{N}(R)$ be
an isometric immersion of a Riemannian product  whose second
fundamental form is adapted to the product net of $M^n$.  Then one of the following possibilities holds:
\begin{itemize}
\item[$(i)$] there exist  an orthogonal decomposition $\Les^N=\Les^{m_0}\times \Pi_{a=1}^r\R^{m_a}$ of Lorentzian space and isometric immersions $f_0\colon\, M_0^{n_0}\to \Hy^{m_0-1}(R_0) \subset \Les^{m_0}$ and $f_a\colon\, M_a^{n_a}\to \Sf^{m_a-1}(R_a)\subset \R^{m_a}$, $1\leq a\leq r$, with 
$-R^2=-R_0^2+\sum_{a=1}^r R_a^2$, such that $f$ is the extrinsic product of
$f_0, \ldots, f_r$. 
\item[$(ii)$] there exist  an orthogonal decomposition
$\R^{N-1}=\Pi_{i=0}^r \R^{m_i}$ and  isometric immersions $f_i\colon M_i\to \R^{m_i}$, $0\leq i\leq r$, such 
that $f$ is the composition $f=j\circ \tilde f$ of $\tilde f=f_0\times\cdots\times f_r
$ with the umbilical inclusion $j\colon\, \R^{N-1}\to \Hy^N(R)$.
\end{itemize}
\end{theorem} 

If we regard each $f_i$, $0\leq i\leq r$,  as an isometric immersion into $\R^{N-1}$ in part $(ii)$ and consider its composition $\tilde f_i=j\circ f_i$ with the umbilical inclusion $j$, then we also say that $f$ is the extrinsic product of $\tilde f_0, \ldots, \tilde f_r$.

\subsection{Partial tubes}

Now we recall a special case of the notion of partial tube introduced in \cite{cw} (see also \cite{cd}).\vspace{1ex}

Let $f_1\colon M_1\to \R^N$ be an  immersion. Assume  there is
an orthonormal set $\{\xi_1, \ldots, \xi_s\}$ of parallel
normal vector fields along $f_1$. In particular, the vector subbundle 
$E=\spa\{\xi_1, \ldots, \xi_s \}$ of
$N_{f_1}M_1$ is parallel and flat, hence the
map $\phi\colon M_1\times \R^s\to E$ given by
$$
\phi_{p_1}(y)=\phi(p_1,y)=\sum_{i=1}^sy_i\xi_i(p_1)
$$ 
for $p_1\in M_1$ and
$y=(y_1, \ldots, y_s)\in \R^s$, is a parallel vector bundle isometry.  Given an 
 immersion $f_0\colon  M_0\to \R^s$, let  $f\colon  M_0\times
M_1\to \R^N$ be defined by
$$
f(p_0, p_1)=f_1(p_1)+\phi_{p_1}(f_0(p_0)).
$$
If  $f$ is an immersion, then it is called the \emph{partial tube over $f_1$ with  fiber $f_0$}, or the \emph{partial tube determined by} $(f_0, f_1, \phi)$.

\begin{remark}\label{re:subst} \po{\em  If  $f$ is a partial tube determined by $(f_0, f_1, \phi)$ as above, then one can always assume that $f_0$ is a \emph{substantial} immersion, i.e., that $f_0(M_0)$ is not contained in any affine subspace of $\R^s$. For if $f_0(M_0)$ is contained 
in the affine subspace $v+\R^\ell\subset \R^s$ then one can replace $f_1$ by its parallel immersion 
$$p_1\in M_1\mapsto f_1(p_1)+\phi_{p_1}(v),$$ and then 
$\phi$  by its restriction to $M_1\times \R^\ell$.}
\end{remark}

\begin{proposition}\po \label{prop:cu} With $f_0,f_1, f$ and $\phi$ as above,  the
following holds:
\begin{itemize}
\item[(i)] The differential of $f$ at $p=(p_0, p_1)$ is given by 
\be\label{eq:diff2}
f_*\tau_0^{p}\!{}_*\,X_0=\phi_{p_1}({f_0}_*X_0),\;\;\;\;\mbox{for}\;X_0\in T_{p_0}M_0,
\ee
and
\be\label{eq:diff1}
f_*\tau_1^{p}\!{}_*\,X_1={f_1}_*(I-A^{f_1}_{\phi_{p_1}(f_0(p_0))})X_1,\;\;\;\;\mbox{for}\;X_1\in T_{p_1}M_1.
\ee 
\item[(ii)] The map $f$  is an immersion at $p$ if and only if
$$
P(p_0,p_1)=I-A^{f_1}_{\phi_{p_1}(f_0(p_0))}
$$ 
is an invertible endomorphism of $T_{p_1}M_1$.
\item[(iii)] If $f$ is an immersion at $p$ then
$$
N_fM(p)=E(p_1)^\perp\oplus 
\phi_{p_1}( N_{f_0}M_0(p_0))\subset N_{f_1}M_1(p_1),
$$
where $E(p_1)^\perp$ is the orthogonal complement of $E(p_1)$ in
$N_{f_1}M_1(p_1)$.
\item[(iv)] If $f$ is an immersion at $p$ then 
\be\label{eq:alphapt}
A^f_\xi(p)\tau_1^{p}\!{}_*=\tau_1^{p}\!{}_*\,P(p_0,p_1)^{-1}A^{f_1}_\xi(p_1),\,\,\,\,\,\mbox{for}\,\,\xi\in N_fM(p),
\ee
\be\label{eq:alphapt2}
A^f_\delta(p)\tau_0^{p}\!{}_*=0, \,\,\,\,\,\mbox{for}\,\,\delta\in E(p_1)^\perp,
\ee
and 
\be\label{eq:alphapt3}
A^f_{\phi_{p_1}(\zeta)}(p)\tau_0^{p}\!{}_*=\tau_0^{p}\!{}_*\,A^{f_0}_{\zeta}(p_0), \,\,\,\,\,\mbox{for}\,\,\zeta\in N_{f_0}M_0(p_0).
\ee
\end{itemize}
\end{proposition} 
\proof The proofs of (\ref{eq:diff2}) and (\ref{eq:diff1}) are straightforward, and the assertions in $(ii)$ and $(iii)$ are immediate 
consequences of those formulas. To prove (\ref{eq:alphapt}), given $\xi\in N_fM(p)$ and $X_1\in T_{p_1}M_1$, let $\gamma\colon\,J\to M_1$ be a smooth curve with $0\in J$, $\gamma(0)=p_1$ and $\gamma'(0)=X_1$. Let $\xi(t)$ be the  parallel transport of $\xi$ along the curve $\tau_1^{p}\circ \gamma$. 
Then, using (\ref{eq:diff1}) we obtain
$$\begin{array}{l}-{f}_*(p)A^{f}_\xi(p) \tau_1^{p}\!{}_*\,X_1=\tilde\nabla_{\tau_1^{p}\!{}_*\,X_1}\xi=\frac{d}{dt}|_{t=0}\xi(p_0,\gamma(t))=
-{f}_1{}_*(p_1)A^{f_1}_\xi(p_1) X_1\vspace{1ex}\\\hspace*{19ex}=-f_*(p_0,p_1)\tau_1^{p}\!{}_*\,P(p_0,p_1)^{-1}A^{f_1}_\xi(p_1) X_1. \end{array}$$
The proofs of (\ref{eq:alphapt2}) and (\ref{eq:alphapt3}) are similar.\vspace{2ex}\qed

 As a consequence of part $(ii)$ of the preceding proposition, we have the following necessary and sufficient condition for $f$ to be an immersion.

\begin{corollary}\po\label{cor:immersion} The map $f$ is an immersion if and only if $f_0(M_0)\subset \Omega(f_1;\phi)$, where 
$$\Omega(f_1;\phi)=\{Y\in \R^s\,:\,(I-A^{f_1}_{\phi_{p_1}(Y)}) \mbox{ is nonsingular for any}\,\,\, p_1\in M_1\}.$$
\end{corollary}

To give a better description of the subset $\Omega(f_1;\phi)$, let  $\eta_1, \ldots, \eta_k\in E$ be the distinct principal normal vector fields of $f_1$ with respect to $E$. Thus, there exists an orthogonal decomposition $TM_1=\oplus_{i=1}^kE_i$ such that $A^{f_1}_\zeta|_{E_i}=\<\zeta, \eta_i\>I$ for any $\zeta\in E$. Therefore, $(I-A^{f_1}_{\phi_{p_1}(Y)})$ is nonsingular if and only if
$\<\phi_{p_1}(Y), \eta_i(p_1)\>\neq 1$ for any $1\leq i\leq k$, that is, if and only if $\phi_{p_1}(Y)$ does not belong to any of the focal hyperplanes 
$$H^{f_1}_i(p_1)=\{\zeta\in N_{f_1}M_1(p_1)\,:\, \<\zeta, \eta_i(p_1)\>=1\}, \,\,\,1\leq i\leq k.$$ 
We conclude that
$$\Omega(f_1;\phi)=\R^s\setminus \left(\bigcup_{p_1\in M_1}\bigcup_{i=1}^k \phi_{p_1}^{-1}(H^{f_1}_i(p_1))\right).$$

Another consequence of Proposition \ref{prop:cu} is the following. 

\begin{corollary}\po\label{cor:pt0} The metric induced by $f$ is a polar metric. More precisely, $$g=\pi_0^*g_0+\pi_1^*(g_1\circ \pi_0),$$
where $g_0$ is the metric on $M_0$ induced by $f_0$ and, for any $p=(p_0, p_1)\in M_0\times M_1$,
the metric $g_1(p_0)$ on $M_1$ is given in terms of the metric $g_1$ induced by $f_1$ by
\be\label{eq:g1p0}
 g_1(p_0)(X_1, Y_1)=g_1(P(p_0,p_1)X_1, P(p_0,p_1)Y_1)
 \ee
for all $X_1, Y_1\in T_{p_1}M_1$. Moreover, the second fundamental form of $f$ is adapted to the product net of $M_0\times M_1$.
\end{corollary}

\subsection[The decomposition theorem]{The decomposition theorem}

The following converse of Corollary \ref{cor:pt0} (and its general version in  Subsection $3.5$ ) is the main result of this paper.

\begin{theorem}\po \label{thm:ptubes0} 
Let $M^n=M_0\times M_1$ be a product manifold endowed with 
a polar metric. 
Let $f\colon M^n\to\R^N$ be an isometric immersion whose second fundamental form
is adapted to the product net. 
Then, there exist an immersion 
$f_1\colon\, M_1\to \R^N$, a parallel vector bundle isometry 
$\phi\colon M_1\times \R^s\to E$ onto a flat parallel subbundle of $N_{f_1}M_1$ and an  
immersion $f_0\colon M_0\to \Omega(f_1; \phi)\subset \R^s$ such
that $f$ is the partial tube determined by $(f_0, f_1, \phi)$.
\end{theorem}

For the proof of Theorem \ref{thm:ptubes0} we need the following lemma.

\begin{lemma}\po \label{prop:start}  Let $f\colon M^n\to\R^N$ be an isometric
immersion and let  $D$ be a vector subbundle of $TM$. Then the 
following conditions on $D$ are equivalent:
\begin{itemize}
\item[(i)] $D$ is totally geodesic and $\alpha_f$ is adapted to $(D,
D^\perp)$.
\item[(ii)] $D$ is integrable and $f_*D^\perp$ is constant in $\R^N$ along each leaf of $D$. 
\end{itemize}
\end{lemma}

\proof If $D$ is integrable, then the subbundle $f_*D^\perp$ is constant in $\R^N$ along each leaf of $D$
if and only if 
$$
\tilde{\nabla}_Xf_*Y\in f_*D^\perp
$$ 
for any $X\in \Gamma(D)$ and $Y\in \Gamma(D^\perp)$. Since 
$$
\tilde{\nabla}_Xf_*Y=f_*\nabla_XY+\alpha(X,Y),
$$
this is the case if and only if $\nabla_XY\in \Gamma(D^\perp)$ and $\alpha_f(X,Y)=0$, that is, if and only if the conditions in  $(i)$ hold.
\vspace{1,5ex}\qed

\noindent \emph{Proof of Theorem \ref{thm:ptubes0}:} 
For a fixed $\bar p_0\in M_0$, define  $f_1\colon M_1\to \R^N$ by  $f_1=f\circ \mu_{\bar p_0}$. 
Let $\E=(E_0, E_1)$ be the product net of $M^n$. 
Given $p_1\in M_1$, it follows from Lemma \ref{prop:start} that the image by $f$
of the leaf $M_0\times \{p_1\}$ of $E_0$ is contained in the affine normal space
of $f_1$ at $p_1$, that is, 
$$
f(p_0,p_1)\in f_1(p_1)+N_{f_1}M_1(p_1)\,\,\,\,\mbox{for every}\,\,p_0\in M_0.
$$
Hence, for each $p_0\in M_0$ we can regard 
 $$p_1\in M_1\mapsto \xi^{p_0}(p_1):=f(p_0,p_1)-f_1(p_1)$$ as
a normal vector field 
along $f_1$. Let $X_1\in T_{p_1}M_1$. Again from Lemma \ref{prop:start},  we
obtain 
$$
\tilde \nabla_{X_1}
\xi^{p_0}=f_*(p_0,p_1)\mu_{p_0}{}_*\,X_1-f_*(\bar{p}_0,p_1)\mu_{\bar p_0}{}_*\,X_1\in f_*(\bar{p}_0,p_1)E_1(\bar{p}_0,p_1)={f_1}_*T_{p_1}M_1.
$$
Hence $\xi^{p_0}$ is a parallel normal vector field along $f_1$. For a fixed $p_1\in
M_1$, set 
$$
E(p_1)=\spa\{\xi^{p_0}(p_1):p_0\in M_0\}.
$$ 
Then, for any pair of points $p_1, \tilde p_1\in M_1$, parallel transport in the normal 
connection of $f_1$ along any curve joining $p_1$ and $\tilde p_1$ takes $E(p_1)$ onto $E(\tilde p_1)$. Thus, such subspaces define a parallel 
flat normal subbundle $E$ of $N_{f_1}M_1$, and hence there exists  a parallel vector bundle isometry $\phi\colon M_1\times \R^s\to E$.

Define $f_0\colon M_0\to \R^s$ by  $\phi_{p_1}(f_0(p_0))=\xi^{p_0}(p_1)$.
Then $$f(p_0, p_1)=f_1(p_1)+\phi_{p_1}(f_0(p_0))$$ for all $(p_0, p_1)\in M^n$. Moreover,  from 
$$f_*(p_0,p_1)\mu_{p_1}{}_*\,X_0=\phi_{p_1}({f_0}_*X_0)\;\;\mbox{for any}\;X_0\in T_{p_0}M_0$$
it follows that $f_0$ is an  immersion.\vspace{2ex}\qed

In view of Theorem \ref{cor:isodrh}, one can also state  Theorem \ref{thm:ptubes0} as follows. 

\begin{corollary}\po \label{cor2:ptubes} 
Let $f\colon\, M^n\to \R^N$ be an isometric immersion of a Riemannian manifold that  carries an orthogonal net 
${\cal E}=(E_0, E_1)$ with $E_0$  totally geodesic.  Assume that the second fundamental form of $f$ is adapted to ${\cal E}$.  Then, there exist locally  a product representation $\psi\colon\; M_0\times M_1\to M$ of $\mathcal{E}$, an immersion $f_1\colon\, M_1\to \R^N$, 
a parallel vector bundle isometry 
$\phi\colon M_1\times \R^s\to E$ onto a flat parallel subbundle of $N_{f_1}M_1$ and an immersion  
 $f_0\colon M_0\to \Omega(f_1; \phi)\subset \R^s$ such
that $f\circ \psi$ is the partial tube determined by $(f_0, f_1, \phi)$. 
\end{corollary}

For instance, Corollary \ref{cor2:ptubes} implies that any surface  with flat normal bundle in $\R^N$ whose curvature lines of one family are geodesics is locally generated by parallel transporting a smooth curve in a normal space of another smooth curve with respect to the normal connection of the latter. Surfaces in $\R^3$ with this property are known in the classical literature as \emph{moulding surfaces}. More precisely, we have the following.

\begin{corollary}\po \label{cor3:ptubes} 
Let $f\colon\, M^2\to \R^N$ be a surface  with flat normal bundle free of umbilical points. Let ${\cal E}=(E_0, E_1)$ be the orthogonal net on $M^2$ determined by its curvature lines. Assume that those correspondent to $E_0$ are geodesics.   Then, there exist locally  a product representation $\psi\colon\; I\times J\to M^2$ of $\mathcal{E}$, where $I,J\subset \R$ are open intervals, a smooth curve $\beta\colon\, J\to \R^N$, 
a parallel vector bundle isometry 
$\phi\colon J\times \R^s\to E$ onto a flat parallel subbundle $E$ of $N_{\beta}J$ and a  smooth curve
 $\alpha\colon I\to \Omega(\beta; \phi)\subset \R^s$ such
that $f\circ \psi$ is the partial tube determined by $(\alpha, \beta, \phi)$. 
\end{corollary}

\subsection[Partial tubes in the sphere and hyperbolic space]{Partial tubes in the sphere and hyperbolic space}

The results of the previous sections can be easily extended to immersions into the sphere and hyperbolic space.

Let $\Ee^{N}$ denote either Euclidean space $\R^{N}$ or Lorentzian space $\Les^{N}$.
Denote by $\Q_\e^{N-1}\subset \Ee^N$  either the sphere  $\Sf^{N-1}$ or hyperbolic space
 $\Hy^{N-1}$, according as $\e=1$ or $\e=-1$, respectively.
 Let $f_1\colon M_1\to \Ee^N$ be an immersion such that 
  $f_1(M_1)$ is contained in $\Q_\e^{N-1}\subset \Ee^N$. 
  Suppose  $E$  is a parallel flat normal subbundle of $N_{f_1}M_1$ with rank $s$ having the position vector field $f_1$ as a section and let $\phi\colon M_1\times \Ee^s\to E$ be a parallel vector bundle isometry. 
   Let $f_0\colon M_0\to \Ee^s$ be an immersion with $$f_0(M_0)\subset (e_1+\Omega(f_1, \phi))\cap \Q_\e^{s-1},$$ where  $e_1\in \Ee^s$ is such that $f_1(p_1)=\phi_{p_1}(e_1)$ for any $p_1\in M_1$. Define  $f\colon M_0\times M_1\to \Ee^N$ by
  \be\label{eq:ptubessphere}f(p_0, p_1)=\phi_{p_1}(f_0(p_0)).\ee

   \begin{corollary}\po\label{cor:ptsphere} The map $f$ is an immersion into $\Q_\e^{N-1}$ whose induced metric is polar  and whose second fundamental form
is adapted to the product net. Conversely, any  immersion $f\colon M_0\times M_1\to\Q_\e^{N-1}\subset \Ee^N$  with these properties is given in this way.
\end{corollary}  
\proof Since $f_0(M_0)\subset \Q_\e^{s-1}$ and $\phi_{p_1}\colon\, \Ee^s\to E(p_1)$ is a linear isometry for any $p_1\in M_1$, it follows that $f(M)\subset  \Q_\e^{N-1}$. We can write
 \be\label{eq:pt*}f(p_0, p_1)=f_1(p_1)+\phi_{p_1}(\tilde f_0(p_0)),\ee
where $\tilde f_0(p_0)=f_0(p_0)-e_1$. Then, the condition $f_0(M_0)\in (e_1+\Omega(f_1, \phi))$ is equivalent to 
$\tilde f_0(M_0)\in \Omega(f_1, \phi)$. Hence $f$ is an immersion, the partial tube determined by $(\tilde f_0, f_1, \phi)$. In particular, the metric induced by $f$ is polar and the second fundamental form of $f$ 
is adapted to the product net of $M_0\times M_1$  by Corollary \ref{cor:pt0}.

 Conversely, if $f\colon M_0\times M_1\to\Q_\e^{N-1}\subset \Ee^N$ is an immersion with these properties, by Theorem \ref{thm:ptubes0} there exist an immersion 
$f_1\colon\, M_1\to \Ee^N$, a parallel vector bundle isometry 
$\phi\colon M_1\times \Ee^s\to E$ onto a flat parallel subbundle of $N_{f_1}M_1$ and an  
immersion $\tilde f_0\colon M_0\to \Omega(f_1; \phi)\subset \Ee^s$ such
that $f$ is the partial tube determined by $(\tilde f_0, f_1, \phi)$. Thus $f$ is given by (\ref{eq:pt*}).

Moreover, by the proof of Theorem~\ref{thm:ptubes0} we may take $f_1$ as
$f_1=f\circ \mu_{\bar p_0}$ for some $\bar p_0\in M_0$, thus $f_1(M_1)\subset \Q_\e^{N-1}$. 
Then, since the position vector is 
a parallel normal vector field, we can assume that $E(p_1)$ contains $f_1(p_1)$ for any $p_1\in M_1$.

Setting $f_1(p_1)=\phi_{p_1}(e_1)$ for a unit vector $e_1\in \Ee^s$, we obtain from (\ref{eq:pt*})
that $f$ is given by (\ref{eq:ptubessphere}) for
  $f_0\colon M_0\to \Ee^s$  defined by $f_0=\tilde f_0+e_1$. Finally, since $f(M)\subset \Q_\e^{N-1}$, it follows from  (\ref{eq:ptubessphere}) that  $f_0(M_0)\subset\Q_\e^{s-1}$.
  \vspace{2ex}\qed 
  
  In particular, from Corollary \ref{cor:ptsphere}  one can write down  the analogues of Corollaries \ref{cor2:ptubes} and 
\ref{cor3:ptubes} for the sphere and hyperbolic space, respectively.

\subsection[A general version of the decomposition theorem]{A general version of the decomposition theorem}

In order to state and prove general versions of 
Corollary \ref{cor:pt0} and Theorem \ref{thm:ptubes0}, consider a product 
$$\tilde f=f_1\times \cdots f_r\colon \tilde M=\Pi_{a=1}^r M_a\to \Pi_{a=1}^r \R^{n_a}=\R^N$$
of immersions $f_a\colon\, M_a\to \R^{n_a}$, $1\leq a\leq r$, such that $N_{\tilde f}\tilde M$ carries a flat parallel subbundle $E$. Let $\phi\colon \tilde M\times \R^s\to E$ be a parallel vector bundle isometry,  let $f_0\colon M_0\to \Omega(\tilde f; \phi)\subset \R^s$ be an immersion and let  $f\colon M= \Pi_{i=0}^r M_i\to \R^N$ 
be the partial tube  determined by $(f_0, \tilde f, \phi)$. 
  
\begin{theorem}\po\label{thm:ptubes} 
The metric induced on $M$ by $f$ is polar 
 and the second fundamental form of $f$ is adapted to the product net. Conversely, any immersion $f\colon M=\Pi_{i=0}^r M_i\to\R^N$ 
 with these properties is given in this way.
\end{theorem}
\proof Regard $M$ as the product $M=M_0\times \tilde M$ and denote by $\tilde \pi\colon\, M\to \tilde M$ the projection. By Corollary \ref{cor:pt0}, the metric induced by $f$ is given by
\be\label{eq:ind1} g=\pi_0^*g_0+\tilde\pi^*(\tilde g\circ \pi_0),\ee
where $g_0$ is the metric on $M_0$ induced by $f_0$ and, for any $p=(p_0, \tilde p)\in M=M_0\times \tilde M$,
the metric $\tilde g(p_0)$ on $\tilde M$ is given in terms of the metric $\tilde g$ induced by $\tilde f$ by
$$
 \tilde g(p_0)(\tilde X, \tilde Y)=\tilde g((I-A^{\tilde f}_{\phi_{\tilde p}(f_0(p_0))})\tilde X, (I-A^{\tilde f}_{\phi_{\tilde p}(f_0(p_0))})\tilde Y)
 $$
for all $\tilde X, \tilde Y\in T_{\tilde p}\tilde M$.

    Since $E$ is a flat parallel subbundle of $N_{\tilde f}\tilde M$, so are its projections $E_a$ onto $N_{f_a}M_a$ for $1\leq a\leq r$.
  Then, there exist  parallel vector bundle isometries $\phi^a\colon  M_a\times \R^{m_a}\to E_a$, $1\leq a\leq r$, such that  $\phi$ is the restriction to $\tilde M\times \R^s$ of the 
parallel vector bundle isometry $\tilde \phi\colon\, \tilde M\times \R^\ell\to \oplus_{a=1}^r E_a$, with  $\R^\ell=\Pi_{a=1}^r \R^{m_a}$,  given by 
\be\label{eq:pvbi}\tilde \phi_{\tilde p}\sum_{a=1}^r v_a=\sum_{a=1}^r \phi^a_{\tilde p_a} v_a,\,\,\,\,\mbox{for any}\,\,\,\,\tilde p=(\tilde p_1, \ldots, \tilde p_r)\in \tilde M.\ee

We denote by  $\tilde \pi_a$ either of the projections $\oplus_{a=1}^r E_a\mapsto E_a$,  $\Pi_{a=1}^r \R^{m_a}\mapsto \R^{m_a}$ and $\tilde M\mapsto M_a$. Also, given $\tilde p=(\tilde{p}_1, \ldots, \tilde{p}_r)\in \tilde M$,  then  ${\tilde\tau}_a^{\tilde{p}}\colon M_a\to \tilde M$ stands for the inclusion given by
 $${\tilde\tau}_a^{\tilde{p}}(p_a)=(\tilde{p}_1, \ldots,  p_a, \ldots, \tilde{p}_r).$$
 From (\ref{eq:so}) we have
$$\tilde f_*A^{\tilde f}_{\phi_{\tilde p}(f_0(p_0))}{\tilde\tau}_a^{\tilde{p}}{}_*={f_a}_*A^{f_a}_{\tilde \pi_a\phi_{\tilde p}(f_0(p_0))}= {f_a}_*A^{f_a}_{\phi^a_{\tilde p_a}(\tilde \pi_a(f_0(p_0)))},\,\,\,1\leq a\leq r. $$
Therefore, 
$$\begin{array}{l}\tilde f_*A^{\tilde f}_{\phi_{\tilde p}(f_0(p_0))}\tilde X=\sum_{a=1}^r\tilde f_*A^{\tilde f}_{\phi_{\tilde p}(f_0(p_0))}\tilde X^a=
\sum_{a=1}^r\tilde f_*A^{\tilde f}_{\phi_{\tilde p}(f_0(p_0))}{\tilde\tau}_a^{\tilde{p}}{}_*{\tilde{\pi}_a}{}_*\tilde X^a\vspace{1.5ex}\\\hspace*{14.2ex}=
\sum_{a=1}^r {f_a}_*A^{f_a}_{\phi^a_{\tilde p_a}(\tilde \pi_a(f_0(p_0)))}{\tilde{\pi}_a}{}_*\tilde X^a.\end{array}$$
We obtain that 
\be\label{eq:ind2}\tilde g(p_0)=\sum_{a=1}^r \tilde \pi_a^* g_a(p_0),\ee
where $g_a(p_0)$ is the metric on $M_a$ given in terms of the metric $g_a$ induced by $f_a$ by 
$$g_a(p_0)(X_a, Y_a)=g_a((I-A^{f_a}_{\phi^a_{\tilde p_a}(\tilde \pi_a(f_0(p_0)))})X_a, (I-A^{f_a}_{\phi^a_{\tilde p_a}(\tilde \pi_a(f_0(p_0)))})Y_a)$$
for all $X_a, Y_a\in T_{\tilde p_a}M_a$.

Since $\tilde \pi_a\circ \tilde \pi=\pi_a$ for $1\leq a\leq r$, we conclude from (\ref{eq:ind1}) and (\ref{eq:ind2}) that
$$g=\pi_0^*g_0+\sum_{a=1}^r\pi_a^*(g_a\circ \pi_0).$$

The assertion  on the second fundamental form of $f$ is a consequence of part $(iv)$ of Propositions \ref{prop:prod} and \ref{prop:cu}, and the proof of the direct statement is completed.

 For the converse, fix $\bar p_0\in M_0$ and  set $\tilde f=f\circ \mu_{\bar
p_0}$. Let $\bar \E=(\bar E_a)_{a=1,\ldots,r}$ be the product net of
$M_{\perp_0}$. Since $\E$ is an orthogonal net and $E_a^\perp$ is totally
geodesic for $1\leq a\leq r$, these properties are inherited by  $\bar \E$, that
is, $\bar \E$ is an orthogonal net and $\bar E_a^\perp$ is totally geodesic for
$1\leq a\leq r$ with respect to the metric on $M_{\perp_0}$ induced by $\mu_{\bar
p_0}$. As  observed after Corollary \ref{cor:prodwarp}, 
this metric  is a product metric.

The second fundamental form $\alpha_{\tilde f}$ of $\tilde f$ is given by
$$
\alpha_{\tilde f}(X,Y)=\alpha_{f}(\mu_{\bar p_0}{}_*\,X,\mu_{\bar p_0}{}_*\,Y)
+f_*\,\alpha_{\mu_{\bar p_0}}(X,Y).
$$
Since  $\<\nabla_{X_b}X_a, X_0\>=0$ for all $X_a\in E_a, X_b\in E_b$, $1\leq
a\neq b\leq r$, and $X_0\in E_0$, because $E_a^\perp$ is totally geodesic for
$1\leq a\leq r$,  we have that $\alpha_{\mu_{\bar p_0}}$ is 
adapted to $\bar \E$. Using this and the fact that
$\alpha_f$ is adapted to $\E$, it follows that
$\alpha_{\tilde f}$ is adapted to $\bar \E$. Hence $\tilde f$ is a 
product of isometric immersions by Moore's Theorem \ref{thm:moore2}.
 
Finally, since $E_0=\cap_{a=1}^r E_a^\perp$ is totally geodesic,   regarding $M^n$ as a product
of $M_0$ and $M_{\perp_0}$ we obtain from  Theorem \ref{thm:ptubes0} that there exist  a parallel vector bundle isometry 
$\phi\colon \tilde M\times \R^s\to E$ onto a flat parallel subbundle of $N_{\tilde f}\tilde M$ and an  
immersion $f_0\colon M_0\to \Omega(\tilde f; \phi)\subset \R^s$ such that $f$ is the partial tube  determined by $(f_0, \tilde f, \phi)$. \vspace{2ex}\qed

As a consequence of Theorem \ref{thm:ptubes}, we classify all local representations of Euclidean space as a  
product manifold endowed with 
a polar  metric. Given an isometric immersion $\tilde f\colon\, M\to \R^N$, we denote by $G\colon\, N_{\tilde f}\tilde M\to \R^N$ the
end-point map given by $$G(v)=\tilde f(\pi(v))+v,$$ where $\pi\colon\,N_{\tilde f}\tilde M\to \tilde M$ is the canonical projection.

\begin{corollary}\po \label{cor:ptubes} 
Let $f\colon M^N=\Pi_{i=0}^r M_i\to\R^N$ be a local isometry of a product manifold endowed with 
a polar  metric. Then, there exist a  product 
$$\tilde f=f_1\times \cdots f_r\colon \tilde M=\Pi_{a=1}^r M_a\to \R^N$$
of immersions  with flat normal bundle, a parallel vector bundle isometry 
$\phi\colon  \tilde M\times \R^s\to N_{\tilde f}\tilde M$  and a local isometry 
$f_0\colon M_0\to \Omega(\tilde f; \phi)\subset \R^s$ such $f=G\circ\phi\circ (f_0\times id)$.
\end{corollary}

It follows from Corollary \ref{cor:ptubes} that any flat polar metric  on a simply connected product manifold $M^N=\Pi_{i=0}^r M_i$ arises as the metric induced (on the open subset of regular points) by the end-point map of the normal bundle of a product of immersions 
$\tilde f=f_1\times \cdots f_r\colon \tilde M=\Pi_{a=1}^r M_a\to \R^N$ with flat normal bundle.\vspace{1ex}

We conclude this subsection by stating the counterparts of Theorem \ref{thm:ptubes} and Corollary~\ref{cor:ptubes} for the sphere and hyperbolic space. Proofs are left to the reader.

 Given an extrinsic product $$\tilde f\colon\, \tilde M=\Pi_{a=1}^r M_a\to \Q_\e^{N-1}\subset \Ee^N,$$ let 
 $\phi\colon \Ee^s\times \tilde M\to E$ be a parallel vector bundle isometry  onto a flat parallel subbundle of 
 $N_{\tilde f}\tilde M$ having the position vector as a section. Let $e_1\in \Ee^s$ be such that $\phi_{\tilde p}(e_1)=\tilde f(\tilde p)$ for any $\tilde p\in \tilde M$, and let $f_0\colon M_0\to \Ee^s$ be an immersion such that $$f_0(M_0)\subset (e_1+\Omega(\tilde f; \phi))\cap \Q_\e^{s-1}.$$ Define $f\colon\, M=\Pi_{i=0}^r M_i=M_0\times \tilde M\to \Ee^N$ by
 \be\label{eq:def}
 f(p_0, \tilde p)=\phi_{\tilde p}(f_0(p_0)).\ee
 
 \begin{theorem}\po\label{thm:ptubeshyp} 
The map $f$ is an immersion taking values in $\Q_\e^{N-1}$  whose  induced metric is polar and whose  second fundamental form  is adapted to the product net of $M$. Conversely, any  map $f\colon \Pi_{i=0}^r M_i\to\Q_\e^{N-1}$ 
 with these properties 
is given in this way.
\end{theorem}

The next corollary implies that any polar metric of constant sectional curvature $\e\in \{-1,1\}$ on an $N$-dimensional  
simply connected product manifold  $\Pi_{i=0}^r M_i$ arises as the induced metric on the (open subset of regular points of the) unit normal 
bundle of an extrinsic  product 
$\tilde f\colon \Pi_{a=1}^r M_a\to \Q_\e^N\subset \Ee^{N+1}$  with flat normal bundle.

\begin{corollary}\po \label{cor:ptubeshyp} 
Let $f\colon M^N=\Pi_{i=0}^r M_i\to\Q_\e^N$ be a local isometry of a product manifold  with 
a polar  metric. 
Then, there exist an extrinsic  product 
$\tilde f\colon \tilde M=\Pi_{a=1}^r M_a\to \Q_\e^N\subset \Ee^{N+1}$  with flat normal bundle, a parallel vector bundle isometry 
$\phi\colon \tilde M\times \Ee^s \to N_{\tilde f}\tilde M$ and a local isometry 
$f_0\colon M_0\to (e_1+\Omega(\tilde f; \phi))\cap \Q_\e^{s-1}\subset \Ee^s$,   where $\phi_{\tilde p}(e_1)=\tilde f(\tilde p)$ for any $\tilde p\in \tilde M$, such that $f$ is given by (\ref{eq:def}).
\end{corollary}

\section[Isometric immersions of quasi-warped products]{Isometric immersions of quasi-warped products}

In this section we use Theorem \ref{thm:ptubes0} to prove a decomposition theorem for isometric immersions of quasi-warped products. Then, we derive  N\"olker's theorem on isometric immersions of warped products follows as a special case.  First we discuss how products and warped products of  immersions can be regarded as  partial tubes.

\subsection{Warped products of immersions}

Let $f_1\colon M_1\to \R^N$ be an isometric immersion such that 
  $f_1(M_1)$ lies  in a subspace $\R^{N-s}\subset \R^N$,  $1\leq s\leq N-1$.  Let $E$ be the vector subbundle of $N_{f_1}M_1$ whose fiber at any $p_1\in M_1$ is the orthogonal complement $\R^{s}$   of $\R^{N-s}$ in $\R^N$, and consider the obvious parallel vector bundle isometry $\phi\colon M_1\times \R^s\to E$. Notice that $\Omega(f_1, \phi)=\R^s$. Given an isometric immersion $f_0\colon M_0\to  \R^s$, the partial tube $f\colon\, M_0\times M_1\to \R^N$  determined by $(f_0, f_1, \phi)$ is just the product $f_0\times f_1\colon\, M_0\times M_1\to \R^s\times \R^{N-s}=\R^N$.

 Now consider an isometric immersion  $f_1\colon (M_1, g_1)\to \R^N$  such that
  $f_1(M_1)$ is contained in   $\Sf^{N-s}=\Sf^{N-1}\cap \R^{N-s+1}$,  $1\leq s\leq N-1$.  Let $E$ be the vector subbundle of $N_{f_1}M_1$ whose fiber at $p_1\in M_1$ is $\spa\{f_1(p_1)\}\oplus \R^{s-1}$, where $\R^{s-1}$ is the orthogonal complement   of $\R^{N-s+1}$ in $\R^N$. Let $\phi\colon M_1\times \R^s\to E$ be a parallel vector bundle isometry  with $\phi_{p_1}(e_1)=f_1(p_1)$ for some unit vector $e_1\in \R^s$ and any $p_1\in M_1$.  Then, for any $Y\in \R^s$ we have $$A^{f_1}_{\phi_{p_1}(Y)}=-\<Y,e_1\>I,$$ where $I$ is the identity endomorphism of $T_{p_1}M_1$. In particular,  
  $$\Omega(f_1, \phi)=\{Y\in \R^s\,:\, \<Y,e_1\>+1\neq 0\}.$$
  Given an isometric immersion $\tilde f_0\colon (M_0, g_0)\to \Omega(f_1, \phi)\subset \R^s$, let $f\colon M_0\times M_1\to \R^N$ be the partial tube determined by $(\tilde f_0, f_1, \phi)$. Then
 $$f(p_0, p_1)=f_1(p_1)+\phi_{p_1}(\tilde f_0(p_0))=\phi_{p_1}(f_0(p_0)),$$
 where $f_0\colon M_0\to \R^s$ is given by $f_0=\tilde f_0+e_1$. Note that the condition $\tilde f_0(M_0)\subset \Omega(f_1, \phi)$ for $f$ to be an immersion reduces to $\<f_0(p_0), e_1\>\neq 0$ for every $p_0\in M_0$.
 
 \begin{proposition}\po\label{prop:n0} The induced metric $g$ on $M_0\times M_1$  is the  warped product of the metrics $g_0$ and $g_1$  with warping function $\rho\colon M_0\to \R_+$ given by $\rho(p_0)=\<f_0(p_0), e_1\>.$
 \end{proposition}
 \proof  By Corollary \ref{cor:pt0}, we have
 $$g=\pi_0^*g_0+\pi_1^*(g_1\circ \pi_0),$$
 with $g_1(p_0)$ given by (\ref{eq:g1p0}) for every $p_0\in M_0$.
 Thus, it suffices to show that 
 \be\label{eq:gp0} g_1(p_0)=\rho^2(p_0)g_1\,\,\,\mbox{for all}\,\,\,p_0\in M_0.\ee
 Since  $$A^{f_1}_{\phi_{p_1}(\tilde f_0(p_0))}=-\<\tilde f_0(p_0),e_1\>I, \,\,\,\mbox{for all}\,\,(p_0,  p_1)\in M_0\times M_1,$$
we obtain that
$$I-A^{f_1}_{\phi_{p_1}(\tilde f_0(p_0))}=\rho(p_0)I,\,\,\,\mbox{for all}\,\,(p_0,  p_1)\in M_0\times M_1,$$
hence (\ref{eq:gp0}) follows from (\ref{eq:g1p0}).\vspace{2ex}\qed

 The map $f$ is called the \emph{warped product} of $f_0$ (or $\tilde f_0$) and $f_1$. If   $f_1\colon \Sf^{N-s}\to \R^{N-s+1}\subset \R^N$, $2\leq s\leq N-1$,  is the canonical umbilical inclusion, then $f$ is a \emph{rotation submanifold}  with $f_0\colon M_0\to \R^s$ as \emph{profile}. On the other hand, for $s=1$ it reduces to the \emph{cone} over  $f_1\colon M_1\to \Sf^{N-1}$. It is convenient to consider also an ordinary product  as a warped  product of immersions.
 
 Note that $f(M_0\times M_1)\subset \Sf^{N-1}$ if and only if $f_0(M_0)\subset \Sf^{s-1}\subset \R^s$. If $f$, $f_0$ and $f_1$ are regarded as immersions into $\Sf^{N-1}$, $\Sf^{s-1}$ and $\Sf^{N-s+1}$, respectively, then $f$ is also said to be the warped product
  of $f_0$ and $f_1$.
  
  Now assume that  
  $f_1\colon M_1\to \Les^N$ is an isometric immersion such that
  $f_1(M_1)$  is contained in the intersection $\Hy^{N-1}\cap V^{N-s+1}$ of hyperbolic space $\Hy^{N-1}$ with an $(n-s+1)$-dimensional
  subspace $V^{N-s+1}$ of $\Les^N$.  Let $E$ be the vector subbundle of $N_{f_1}M_1$ whose fiber at $p_1\in M_1$ is the Lorentzian subspace $\spa\{f_1(p_1)\}\oplus V^\perp$,  and let $\phi\colon M_1\times \Les^s\to E$ be a parallel vector bundle isometry  with $\phi_{p_1}(e_1)=f_1(p_1)$ for some unit vector $e_1\in \Les^s$ and any $p_1\in M_1$. Let $f_0\colon\, M_0\to \Hy^{s-1}\subset \Les^s$ be an isometric immersion such that $\<f_0, e_1\>\neq 0$ for all $p_0\in M_0$. Then, the map $f\colon M_0\times M_1\to  \Les^N$ given by 
 $$f(p_0, p_1)=\phi_{p_1}(f_0(p_0))$$
 is an immersion and satisfies $f(M_0\times M_1)\subset \Hy^{N-1}$. As in the spherical case, it is called  the warped product
  of $f_0$ and $f_1$. If $f_1$ is the umbilical inclusion of $\Hy^{N-1}\cap V^{N-s+1}$ into $\Hy^{N-1}$, then $f$ is the \emph{rotational 
  submanifold with $f_0$ as profile}. It is said to be of \emph{spherical}, \emph{hyperbolic} or \emph{parabolic} type, according as the
  subspace $V$ is space-like, time-like or degenerate, respectively. Again, extrinsic products of immersions into the sphere and hyperbolic space are also considered as warped products of immersions.
  
\subsection{Partial tubes over curves}

Consider  a unit speed curve $\gamma\colon I\to \R^N$. Let $\phi\colon I\times \R^s\to E$ be a parallel vector bundle isometry onto a parallel flat vector subbundle $E$ of $N_{\gamma}I$.  Since
 $$A^{\gamma}_{\phi_{t}(Y)}=\<\gamma''(t), \phi_{t}(Y)\>I $$
 for all $(t, Y) \in I\times \R^s$, it follows that 
 $$\Omega(\gamma, \phi)=\{Y\in \R^s\,:\, \<\gamma''(t), \phi_{t}(Y)\>\neq 1\,\,\,\,\mbox{for all}\,\,t\in I\}.$$
 Equivalently, denoting by 
 $$H_{t}^{\gamma}=\{Y\in \R^N\,:\,\<Y, \gamma''(t)\>=1\}$$ 
  the focal hyperplane of $\gamma$ at $t$, we have
    $$\Omega(\gamma, \phi)=\R^s\setminus\bigcup_{t\in I}\phi_{t}^{-1}(H^{\gamma}_{t}).$$
  Given an isometric immersion $ f_0\colon (M_0, g_0)\to \Omega(\gamma, \phi)\subset \R^s$, let  
  $f\colon M_0\times M_1\to \R^N$ be the partial tube determined by $( f_0, \gamma, \phi)$. Arguing as in the proof of Proposition \ref{prop:n0}, we obtain the following.

 \begin{proposition}\po\label{prop:n0b} The induced metric on $M_0\times I$  is the quasi-warped product of the metrics $g_0$ and the standard metric on  $I$ with warping function $\rho\colon M_0\times I\to \R_+$ given by $$\rho(p_0,t)=1-\<\gamma''(t), \phi_{t}(f_0(p_0))\>.$$
 \end{proposition}
 
 Now suppose that  $\gamma\colon I\to \Ee^N$ is a unit speed curve taking values in $\Q_\e^{N-1}$.   
  Let $E$  be a parallel flat normal subbundle of $N_{\gamma}I$ of rank $s$ having the position vector $\gamma$ as a section and let $\phi\colon\, I\times \Ee^s\to E$ be a parallel vector bundle isometry such that 
   $\gamma(t)=\phi_{t}(e_1)$ for every $t\in I$ and some unit time-like vector $e_1\in \Ee^s$. Let $f_0\colon (M_0, g_0)\to \Q_\e^{s-1}\subset \Ee^s$ be an isometric immersion such that
   $$\rho(p_0,t):=\<\gamma''(t), \phi_t(f_0(p_0))\>\neq 0$$
   for all $t\in I$ and $p_0\in M_0$.  Define  $f\colon M_0\times I\to \Ee^N$ by 
   $$f(p_0, t)=\phi_t(f_0(p_0)).$$
   We also call $f$ a partial tube over $\gamma$ with fiber $f_0$. 

  \begin{proposition}\po\label{prop:n0c} The map $f$ takes values in $\Q_\e^{N-1}$ and its induced metric   is the quasi-warped product of the metrics $g_0$ and the standard metric on  $I$ with warping function~$\rho$.
 \end{proposition}

\subsection{Decomposition of isometric immersions of quasi-warped products}

  The next result  shows that  the  special cases of partial tubes in the two previous subsections comprise all possible examples of isometric immersions into Euclidean space of a quasi-warped product manifold with  two factors whose second fundamental forms are adapted to the product structure. The case of arbitrarily many factors is considered in the next subsection. In the sequel, if $\e=0$  then $\Q_\e^N$ stands for Euclidean space $\R^N$.
 
 We point out that, if $g=\pi_0^*g_0+\rho^2\pi_1^*g_1$ is a quasi-warped product metric on  $M_0\times M_1$, then we can assume that either  $g$ is a Riemannian product metric or $\rho$ does not depend only on $M_1$. For if 
 $\rho=\bar \rho\circ \pi_1$ for some $\bar \rho\in C^\infty (M_1)$, then  $g$ is the Riemannian product of $g_0$ and the metric  $\bar\rho^2g_1$ on $M_1$. We make this assumption  in the following. 
 
  \begin{theorem}\po\label{prop:n1} Let $f\colon M=M_0\times_{\rho} M_1\to \Q_\e^N$, $\e\in \{-1, 0, 1\}$,  be an isometric immersion of a quasi-warped product  whose   second fundamental form  is adapted to the product net. Then  $f$ is either a warped product of immersions or  a partial tube over a curve.
 \end{theorem}
 \proof  We give the proof for Euclidean space, the others being similar.  If $M$ is a Riemannian product, we conclude from Moore's  theorem that $f$ is a  product of  immersions. In view of the remark in the beginning of this subsection, we assume from now on that the warping function $\rho$ does not depend only on $M_1$.
 
 Fixed $\bar p_0\in M_0$, let  $f_1\colon M_1\to \R^N$ be given by  $f_1=f\circ \mu_{\bar p_0}$. 
 Notice that the metric induced by $f_1$ is $g_1(\bar p_0)=\rho^2_{\bar p_0}g_1,$ where $g_1$ is the metric on $M_1$ and $\rho_{\bar p_0}=\rho\circ \mu_{\bar p_0}$.
  By Theorem \ref{thm:ptubes0}, there exist   a parallel vector bundle isometry $\phi\colon M_1\times \R^s\to E$ onto a flat parallel   subbundle of $N_{f_1}M_1$ and  an   isometric immersion $f_0\colon M_0\to \Omega(f_1; \phi)\subset\R^s$ such that  $f$  is the partial tube determined by $(f_0, f_1, \phi)$.  By Remark \ref{re:subst}, we can assume that $f_0$ is substantial in $\R^s$.
  
  By Corollary \ref{cor:pt0}, the  metric induced by $f$ is given by
$$g=\pi_0^*g_0+\pi_1^*(g_1\circ \pi_0),$$
where $g_0$ is the metric on $M_0$  and, for any $p=(p_0, p_1)\in M_0\times M_1$,
the metric $g_1(p_0)$ on $M_1$ is given in terms of the metric $g_1(\bar p_0)$ induced by $f_1$  by
$$ g_1(p_0)(X_1, Y_1)=g_1(\bar p_0)((I-A^{f_1}_{\phi_{p_1}(f_0(p_0))})X_1, (I-A^{f_1}_{\phi_{p_1}(f_0(p_0))})Y_1)
$$
for all $X_1, Y_1\in T_{p_1}M_1$. Therefore, we must have
$$\rho^2_{p_0}g_1(X_1, Y_1)=g_1(p_0)(X_1, Y_1)=\rho^2_{\bar p_0}g_1((I-A^{f_1}_{\phi_{p_1}(f_0(p_0))})X_1, (I-A^{f_1}_{\phi_{p_1}(f_0(p_0))})Y_1)$$
for all $X_1, Y_1\in T_{p_1}M_1$, hence
  \be\label{eq:fund}(I-A^{f_1}_{\phi_{p_1}(f_0(p_0))})^2=\frac{\rho^2_{p_0}(p_1)}{\rho^2_{\bar p_0}(p_1)}I,\,\,\,\mbox{for all}\,\,\,(p_0, p_1)\in M_0\times M_1.\ee

 It suffices to prove that $E$ is an umbilical subbundle of $N_{f_1}M_1$ and that $E\not\subset N_1^\perp(f_1)$. For this implies that either $M_1$ is one-dimensional or   $f_1(M_1)$ is contained in an $(N-s)$-dimensional sphere $\Sf^{N-s}$, $1\leq s\leq N-1$, which we can assume to be of unit radius and centered at the origin of a subspace $\R^{N-s+1}\subset \R^N$. Moreover, in the latter case  $E$ must be  the vector subbundle  of $N_{f_1}M_1$ whose fiber at $p_1\in M_1$ is spanned by the position vector $f_1(p_1)$ and the orthogonal complement  $\R^{s-1}$ of $\R^{N-s+1}$ in $\R^N$.

    By the assumption on  $\rho$, equation (\ref{eq:fund})   rules out the possibility that $E\subset N_1^\perp(f_1)$. 
    Let $\eta_1, \ldots, \eta_k\in E$ be the distinct principal normal vector fields of $f_1$ with respect to $E$. Thus, there exists an orthogonal decomposition $TM_1=\oplus_{i=1}^kE_i$ such that $$A^{f_1}_\zeta|_{E_i}=\<\zeta, \eta_i\>I$$ for any $\zeta\in E$. We must show that $k=1$. 
    
    Set $\lambda(p_0,p_1)=(\rho^2_{p_0}(p_1)/\rho^2_{\bar p_0}(p_1))^{1/2}$ and write $\eta_i(p_1)={\phi_{p_1}(V_i(p_1))}$. Then (\ref{eq:fund}) can be written as 
    $$(1-\<f_0(p_0),V_i(p_1)\>)^2=\lambda^2(p_0, p_1)    ,\,\,\,\mbox{for all}\,\,\,(p_0, p_1)\in M_0\times M_1\,\,\mbox{and}\,\,1\leq i\leq k.$$
 Thus, for each $1\leq i\leq k$, either $\<V_i(p_1),f_0(p_0)\>=1+\lambda(p_0, p_1)$ or  $\<V_i(p_1),f_0(p_0)\>=1-\lambda(p_0, p_1)$  for all $(p_0, p_1)\in M_0\times M_1$. If $k>2$, then there exist $1\leq i\neq j\leq k$ such that $\<V_i(p_1)-V_j(p_1), f_0(p_0)\>=0$ for any $(p_0, p_1)\in M_0\times M_1$, contradicting the fact that $f_0$ is substantial in $\R^s$. If $k=2$, from  $\<V_1(p_1),f_0(p_0)\>=1+\lambda(p_0,p_1)$ and  $\<V_2(p_1),f_0(p_0)\>=1-\lambda(p_0,p_1)$ we obtain   that $\<f_0(p_0), V_1(p_1)+V_2(p_1)\>=2$ for all $(p_0, p_1)\in M_0\times M_1$, also contradicting the fact that $f_0$ is substantial in $\R^s$. Hence $k=1$.\vspace{1ex}\qed
 
\begin{remark}\po \emph{  In part $(ii)$  we have, in addition, that $f_0$ is an isometric immersion and that there exist $\rho_0\in C^{\infty}(M_0)$ and $\rho_1\in C^{\infty}(M_1)$ such that $\rho=(\rho_0\circ \pi_0)(\rho_1\circ \pi_1)$ and $f_1$ is an isometric immersion with respect to $\tilde g_1=\rho_1^2g_1$.}
\end{remark}

 N\"olker's theorem for isometric immersions into space forms of a warped product manifold with only two factors now follows easily from Theorem \ref{prop:n1}.

 \begin{theorem}\po\label{prop:n1b} Any isometric immersion $f\colon M_0\times_{\rho} M_1\to \Q_\e^N$  of a warped product manifold
 whose second fundamental form is adapted to the product net of $M_0\times M_1$ is  a warped product of immersions. 
 \end{theorem}
 \proof   Again, we give the proof for Euclidean space, the others being similar. It suffices to prove that if $f$ is a partial tube over a curve $\gamma\colon\,I\to \R^N$,  then   $\gamma(I)$ is contained in a sphere of $\R^N$, which we can assume to be a hypersphere  $\Sf^{N-s}$ of  unit radius  centered at the origin of a subspace $\R^{N-s+1}\subset \R^N$, $1\leq s\leq N-1$, and the fiber at $t\in I=M_1$ of the vector subbundle $E$ of $N_{\gamma}I$  is spanned by the position vector $\gamma(t)$ and the orthogonal complement  $\R^{s-1}$ of $\R^{N-s+1}$ in $\R^N$.
  
  We can assume that  $\gamma=f\circ \mu_{\bar p_0}$ for some $\bar p_0\in M_0$, and hence the fact that $\gamma$ has unit speed means that $\rho(\bar p_0)=1$. Then, the assumption on $\rho$ becomes 
  that $\rho$ is not identically one on $M_0$. By Proposition \ref{prop:n0b}, the  warping function $\rho\colon M\to \R_+$ is given by 
  \be\label{eq:rho}
  \rho(p_0,t)=1-\<\gamma''(t), \phi_{t}(f_0(p_0))\>.
  \ee
  Since $\rho$ does not depend on $t$ by assumption,  differentiating $\rho$ with respect to $t$ and using that $f_0$ is substantial yields
  $$\<\gamma'''(t), \xi\>=0\,\,\,\,\mbox{for all}\,\,\,t\in I\,\,\,\mbox{and for all}\,\,\,\,\xi  \in E(t).$$
  Observe also that  one can not have $\gamma''(t)\in E^\perp (t)$ at any $t\in I$, for this and (\ref{eq:rho}) would imply 
  $\rho$ to be identically one.
  
  The result is now a consequence of the next lemma.
  
  \begin{lemma}\po  Let $\gamma\colon\, I\to \R^N$ be a unit-speed curve. Assume that there exists a  parallel normal subbundle $E$ of $N_{\gamma}I$ of rank $s$ such that $\gamma'''(t) \in E(t)^\perp$  but $\gamma''(t)\not \in E(t)^\perp$ for any $t\in I$.   
Then, up to a rigid motion of $\R^N$,  $\gamma(I)$ is contained in a hypersphere  $\Sf^{N-s}(R)$ of  radius $R$  centered at the origin of a subspace $\R^{N-s+1}\subset \R^N$, $1\leq s\leq N-1$, and $E(t)$ is spanned by $\gamma(t)$ and the orthogonal complement  $\R^{s-1}$ of $\R^{N-s+1}$ in $\R^N$.
  \end{lemma}
 \proof  Since $\gamma''(t)\not \in E(t)^\perp$ for any $t\in I$,   the orthogonal projection $(\gamma''(t))_{E(t)}$ of $\gamma''(t)$ onto $E(t)$ is nowhere vanishing. Let $\zeta(t)$ be a unit vector field along $\gamma$ in the direction of
 $(\gamma''(t))_{E(t)}$.   For any section $\xi$ of the orthogonal complement $F=\{\zeta\}^\perp$ of $\{\zeta\}$ in $E$, using that $\gamma'''(t) \in E(t)^\perp$ and that $E$ is parallel in the normal connection of $\gamma$  we obtain
  $$\<\xi', \zeta\>=\<\xi',  \gamma''\>=-\<\gamma''', \xi\>=0.$$ 
  It follows that $F$ is also parallel in the normal connection of $\gamma$, and hence $F$ is  a constant subspace $\R^{s-1}$ of $\R^N$. Hence, up to a rigid motion of $\R^N$, we may assume that $\gamma(I)$ is contained in the  orthogonal complement  $\R^{N-s+1}$ of $\R^{s-1}$ in $\R^N$. Moreover, we have $\zeta'=\lambda \gamma'$, with  $\lambda=\<\zeta', \gamma'\>=-\<\zeta, \gamma''\>$. 
  Now,  $\<\zeta, \gamma''\>'=\<\zeta', \gamma''\>+\<\zeta, \gamma'''\>=0$, 
  hence $\lambda$ is a nonzero constant $1/R\in \R$. It follows that  $\gamma-{R}\zeta$ is a constant vector of $\R^{N-s+1}$, which we can assume to be $0$.\qed

  \subsection[A general version of Theorem \ref{prop:n1}]{A general version of Theorem \ref{prop:n1}}
  
  Let $\tilde M=\Pi_{a=1}^r M_a$ be a product manifold, and let $f_a\colon\, M_a\to \R^{N_a}$, $1\leq a\leq r$, be isometric   immersions.
  Suppose that there exist $1\leq k\leq \ell\leq r$ such that $f_a(M_a)$ is  contained in a hypersphere $\Sf^{N_a-1}(R_a)$ of radius $R_a$ centered at the origin for $1\leq a\leq k$,    and $M_a$ is one-dimensional for $k+1\leq a\leq \ell$.  Given $1\leq a\leq k$ and   $p_a\in M_a$, let $E_a(p_a)$ 
  denote the one-dimensional subspace spanned by the position vector $f_a(p_a)$. For  $k+1\leq a\leq \ell$,  let $E_a$ be any flat parallel subbundle of $N_{f_a}M_a$.
 
  Define $$\tilde f=f_1\times \cdots \times f_r\colon\, \Pi_{a=1}^r M_a\to \Pi_{a=1}^r \R^{N_a}=\R^{N-m}\subset \R^m\times \R^{N-m}=\R^{N},$$ and let $E$ be the flat parallel subbundle of the normal bundle of $\tilde f$ in $\R^{N}$ whose fiber at $p=(p_1,\ldots, p_r)$ is $(\oplus_{a=1}^\ell E_a(p_a))\oplus \R^m$. Choose parallel vector bundle isometries $\phi^a\colon\, M_a\times \R^{n_a}\to V_a$, $k+1\leq a\leq \ell$, set $\R^s=\R^k\times \Pi_{a=k+1}^\ell \R^{n_a}\times \R^m$, choose an orthonormal basis $e_1, \ldots, e_k$ of $\R^k$ and denote by $\pi_a\colon\, \R^s\to\R^{n_a}$ and $\bar\pi\colon\,\R^s\to \R^m$ the projections onto $\R^{n_a}$ and  $\R^m$, respectively. 
  
  Define a parallel vector bundle isometry $\phi\colon \tilde M\times \R^s\to E$ by requiring that $\phi_{\tilde p}(e_a)=f_a(p_a)/R_a$ for  $\tilde{p}=(p_1, \ldots, p_r)$ and $1\leq a\leq k$, 
  $\phi_{\tilde p}|_{\R^{n_a}}=\phi_{p_a}^a$ for $k+1\leq a\leq \ell$ and $\phi_{\tilde p}|_{\R^{m}}$ be any linear isometry of $\R^m$.
  
%
 Notice that $Y\in \sum_{a=1}^kR_ae_a+\Omega(\tilde f; \phi)$ if and only if
  $$\<Y, e_a\>\neq 0,\,\,1\leq a\leq k \,\,\mbox{and}\,\, \pi_a(Y)\not \in (\phi^a_{p_a})^{-1}(\cup_{p_a\in M_a} H^{f_a}_{p_a}),\,\,k+1\leq a\leq \ell.$$
  
  Given an    immersion  $f_0\colon M_0\to (\sum_{a=1}^kR_ae_a+\Omega(\tilde f; \phi))\subset \R^s$, it follows that the immersion $\tilde f_0\colon\, M_0\to \R^s$  defined by $\tilde f_0=f_0-\sum_{a=1}^k R_ae_a$ satisfies$\tilde f_0(M_0)\subset\Omega(\tilde f; \phi)$. Let $f\colon M^n=\Pi_{i=0}^r M_i\to \R^N$ be the partial tube determined by $(\tilde f_0, \tilde f, \phi)$. In case $k=\ell$, the immersion $f$ is called the \emph{warped product} of $f_0, \ldots f_r$.

 \begin{theorem}\po\label{prop:n1c} The metric $g$ induced by $f$ is  the quasi-warped product of the metrics $g_0, \ldots, g_r$ induced by $f_0, \ldots, f_r$, respectively,   with  warping functions $\rho_a\colon M_0\times M_a\to \R_+$ given by
 $$\rho_a(p_0, p_a)=\rho_a(p_0)=\frac{1}{R_a}\<f_0(p_0), e_a\>\,\,\,\mbox{ for}\,\,\,1\leq a\leq k,$$
 $$\rho_a(p_0, p_a)=1-\<f_a''(p_a), \phi_{p_a}(f_0(p_0))\>\,\,\,\mbox{ for}\,\,\,k+1\leq a\leq \ell\,\,\,\,\mbox{and}\,\,\,\rho_a=1\,\,\,\mbox{ for}\,\,\,\ell+1\leq a\leq r.$$
 Moreover,  the second
  fundamental form of $f$ is adapted to the product net of $M^n$. 
 Conversely, any isometric immersion $f\colon M^n\to \R^N$  of a quasi-warped product manifold whose second
  fundamental form is adapted to its product net   is given in this way. 
 \end{theorem}
 \proof  As shown in the proof of Theorem \ref{thm:ptubes}, the metric $g$ induced by $f$ is given by 
 $$g=\pi_0^*g_0+\sum_{a=1}^r\pi_a^*(g_a\circ \pi_0),$$
 where $g_a(p_0)$ is the metric on $M_a$ given in terms of the metric $g_a$ induced by $f_a$ by 
$$g_a(p_0)(X_a, Y_a)=g_a((I-A^{f_a}_{\phi^a_{p_a}(\pi_a(\tilde f_0(p_0)))})X_a, (I-A^{f_a}_{\phi^a_{p_a}(\pi_a(\tilde f_0(p_0)))})Y_a)$$
for all $X_a, Y_a\in T_{p_a}M_a$.

For $1\leq a\leq k$, we have $\pi_a(\tilde f_0(p_0)))=\<\tilde f_0(p_0), e_a\>e_a$, hence
$$A^{f_a}_{\phi^a_{p_a}(\pi_a(\tilde f_0(p_0)))}=-\frac{1}{R_a}\<\tilde f_0(p_0), e_a\>I.$$
Therefore,
$$g_a(p_0)=(1+\frac{1}{R_a}\<\tilde f_0(p_0), e_a\>)^2g_a=\rho^2_a(p_0)g_a.$$
For $k+1\leq a\leq \ell$ we have 
$$A^{f_a}_{\phi^a_{p_a}(\pi_a(\tilde f_0(p_0)))}=\<f_a''(p_a), \phi^a_{p_a}(\tilde f_0(p_0))\>I=\<f_a''(p_a), \phi^a_{p_a}( f_0(p_0))\>I,$$
thus
$$g_a(p_0)=(1-\<f_a''(p_a), \phi^a_{p_a}( f_0(p_0))\>)^2g_a=\rho_{p_0}^2g_a.$$

Finally, since $A^{f_a}_{\phi^a_{ p_a}(\pi_a(\tilde f_0(p_0)))}=0$ for $\ell+1\leq a\leq r$,  for such values of $a$ 
we have $g_a(p_0)=g_a$ for any $p_0\in M_0$. The assertion on the second fundamental form of $f$ follows from Theorem \ref{thm:ptubes}. This completes the proof of the direct statement.

    Conversely, let $f\colon M^n=M_0\times_{\rho} \Pi_{a=1}^r M_a\to \R^N$ be an isometric immersion  of a quasi-warped product manifold whose second fundamental form is adapted to its product net. Thus,  the metric $g$ on $M^n$ has the form
    $$g=\pi_0^*g_0+\sum_{a=1}^r(\rho_a\circ \pi_{0,a})^2\pi_a^*g_a$$
for some $\rho_a\in C^\infty(M_0\times M_a)$ with $\rho_a>0$, $1\leq a\leq r$.

For a fixed $\bar p_0\in M_0$, let  $\tilde f\colon \tilde M=\Pi_{a=1}^r M_a\to \R^N$ be defined by  $\tilde f=f\circ \mu_{\bar p_0}$. Then, the metric induced by $\tilde f$ is 
$$\mu_{\bar p_0}\!{}^*g=\sum_{a=1}^r(\rho_a)_{\bar p_0}^2\tilde \pi_a^*g_a,$$
where $\tilde \pi_a\colon\, \tilde M\to M_a$ is the projection.
Hence, we can replace each $g_a$ by $(\rho_a)_{\bar p_0}^2g_a$ and $\rho_a$ by $\tilde \rho_a$, given by $(\tilde \rho_a)_{p_0}=(\rho_a)_{p_0}/(\rho_a)_{\bar p_0}$ for any $p_0\in M_0$, so as to make $\tilde f$ into an isometric immersion with respect to the product metric of $g_1, \ldots, g_r$ on $\tilde M$.  

Since the second fundamental form of $f$ 
is adapted to the product net of $M^n$, it follows from  Theorem \ref{thm:ptubes} that  $\tilde f\colon \tilde M\to \R^N$ is a product 
 $$\tilde f=f_1\times \cdots \times f_r\colon\, \Pi_{a=1}^r M_a\to \Pi_{a=1}^r \R^{N_a}=\R^{N-\ell}\subset \R^\ell\times \R^{N-\ell}=\R^{N}$$
  of substantial  isometric immersions $f_a\colon M_a\to \R^{N_a}$, and there exist a flat parallel subbundle $E$ of $N_{\tilde f}\tilde M$,  a parallel vector bundle isometry $\phi\colon \tilde M\times \R^s\to E$,  and  an  isometric immersion $\tilde f_0\colon M_0\to \Omega(\tilde f; \phi)\subset\R^s$, which we can assume to be substantial by Remark~\ref{re:subst},  such that  $f\colon M^n\to \R^N$  is the partial tube determined by $(\tilde f_0, \tilde f, \phi)$.

As in the proof of Theorem \ref{thm:ptubes}, let $E_a$ be the projection of  $E$  onto $N_{f_a}M_a$ for $1\leq a\leq r$ and let $\phi^a\colon  M_a\times \R^{m_a}\to E_a$ be a  parallel vector bundle isometry, so that $\phi$ is the restriction of the 
parallel vector bundle isometry $\tilde \phi\colon\, \R^\ell=\Pi_{a=1}^r \R^{m_a}\to \Pi_{a=1}^r E_a$ given by (\ref{eq:pvbi}). 
By Theorem \ref{thm:ptubes}, we have 
  $$g=\pi_0^*g_0+\sum_{a=1}^r\pi_a^*(g_a\circ \pi_0),$$
  where, for any $(p_0, \ldots, p_r)\in M$, 
  $$g_a(p_0)(X_a, Y_a)=g_a((I-A^{f_a}_{\phi^a_{ p_a}(\pi_a(f_0(p_0)))})X_a, (I-A^{f_a}_{\phi^a_{ p_a}(\pi_a(f_0(p_0)))})Y_a)$$
for all $X_a, Y_a\in T_{p_a}M_a$.
    Therefore,  we must have
    $$(I-A^{f_a}_{\phi^a_{p_a}(\pi_a(f_0(p_0)))})^2=(\tilde \rho_a^2)(p_a)I,$$
    and we conclude as in the proof of Theorem \ref{prop:n1} that each $E_a$ either belongs to $N_1^\perp(f_a)$ or it is an umbilical
    subbundle of $N_{f_a}M_a$. The conclusion follows.\vspace{2ex} \qed

 A special case of  Theorem \ref{prop:n1c} is a classification of all local representations of Euclidean space as a  
quasi-warped product, i.e., local isometries of quasi-warped product manifolds into Euclidean space. As before, given an isometric immersion $f\colon\, M^n\to \R^N$ we denote by $G\colon\, N_fM\to \R^N$ the
end-point map given by $G(v)=f(\pi(v))+v$, where $\pi\colon\,N_fM\to M$ is the canonical projection.

\begin{corollary}\po \label{cor:ptubes2} 
Let $\psi\colon M^N=\Pi_{i=0}^r M_i\to\R^N$ be a local isometry of a quasi-warped product manifold. 
Then, there exist a  product 
$$\tilde f=f_1\times \cdots \times f_r\colon\, \tilde M:=\Pi_{a=1}^r M_a\to \Pi_{a=1}^r \R^{N_a}=\R^{N-m}\subset \R^m\times \R^{N-m}=\R^{N},$$
where $f_a\colon\, M_a\to\Sf^{N_a-1}(R_a)\subset \R^{N_a}$ is a local isometry  for $1\leq a\leq k$   and a unit speed curve  for $k+1\leq a\leq  r$,   a parallel vector bundle isometry 
$\phi\colon \R^s\times \tilde M\to N_{\tilde f}\tilde M$  and a local isometry 
$f_0\colon M_0\to  (\sum_{a=1}^kR_ae_a+\Omega(\tilde f; \phi))\subset \R^s$ such that $\psi=G\circ\phi\circ (f_0\times id)$, with $M^N$ regarded as the product $M^N=M_0\times \tilde M$.
\end{corollary}

We call $\psi$ the {\em quasi-warped representation\/} of $\R^n$ determined by $\tilde f$. 
N\"olker's warped product representations correspond to the particular cases in which all  unit-speed curves $f_a$, $k+1\leq a\leq \ell$, are circles. Notice that, in this case, we have
$$ \sum_{a=1}^kR_ae_a+\Omega(\tilde f; \phi)=\Omega(\tilde f):=\{Y\in \R^s:\<Y, e_a\>\neq 0,\,\,1\leq a\leq \ell\}.$$

To describe all isometric immersions $f\colon\, M^n\to \Q_\epsilon^{N-1}$, $\epsilon \in \{-1, 1\}$, of a quasi-warped product manifold whose
second fundamental forms are adapted to the product net of $M^n$, start with an extrinsic product 
$$\tilde f\colon\, \tilde M=\Pi_{a=1}^r M_a\to \Q_\epsilon^m\subset \Q_\e^{N-1}\subset \Ee^N,$$
with $\tilde f$ substantial in $\Q_\epsilon^m$. 
Suppose that there exists $1\leq k\leq  r$ such that  $M_a$ is one-dimensional for $k+1\leq a\leq r$.  
Given $1\leq a\leq k$ and   $p_a\in M_a$, let $E_a(p_a)$ 
  denote the one-dimensional subspace spanned by the position vector $f_a(p_a)$. For  $k+1\leq a\leq r$,
    let $E_a$ be any flat parallel subbundle of $N_{f_a}M_a$ having the position vector $f_a$ as a section, and
let $E_0$ be the normal space of $\Q_\epsilon^m$ in $\Q_\e^{N-1}$. Then $E=\oplus_{i=0}^r E_a$ is a flat parallel subbundle of 
 $N_{\tilde f}\tilde M$ having the position vector $\tilde f$ as a section.
Let 
 $\phi\colon \Ee^s\times \tilde M\to E$ be a parallel vector bundle isometry and let $e_1\in \Ee^s$ be such that $\phi_{\tilde p}(e_1)=\tilde f(\tilde p)$ for any $\tilde p\in \tilde M$. Given an immersion $f_0\colon M_0\to \Ee^s$  such that $$f_0(M_0)\subset (e_1+\Omega(\tilde f; \phi))\cap \Q_\e^{s-1},$$
 define $f\colon\, M=\Pi_{i=0}^r M_i=M_0\times \tilde M\to \Ee^N$ by
 \be\label{eq:def2}
 f(p_0, \tilde p)=\phi_{\tilde p}(f_0(p_0)).\ee
 
 Then we have the following version of Theorem \ref{prop:n1c}  for  the sphere and hyperbolic space.
 
 \begin{theorem}\po\label{prop:n1e} 
The map $f$ is an immersion taking values in $\Q_\e^{N-1}$  whose  induced metric is a quasi-warped product metric and whose  second fundamental form  is adapted to the product net of $M$. Conversely, any  immersion $f\colon \Pi_{i=0}^r M_i\to\Q_\e^{N-1}$ 
 with these properties is given in this way.
\end{theorem}

The counterpart of  Corollary \ref{cor:ptubes2} is as follows. It implies that any quasi-warped product metric of constant sectional curvature $\e\in \{-1,1\}$ on an $N$-dimensional  
simply connected product manifold  $\Pi_{i=0}^r M_i$ arises as the induced metric on the (open subset of regular points of the) unit normal 
bundle of an extrinsic  product 
$$\tilde f=f_1\times \cdots\times f_r\colon \Pi_{a=1}^r M_a\to \Q_\e^N\subset \Ee^{N+1},$$ with
$f_a$ spherical or $M_a$ one-dimensional for every $1\leq a\leq r$.

\begin{corollary}\po \label{cor:ptubeshyp2} 
Let $f\colon M^N=\Pi_{i=0}^r M_i\to\Q_\e^N$ be a local isometry of a quasi-warped product manifold. 
Then, there exist an extrinsic  product 
$$\tilde f=f_1\times \cdots\times f_r\colon \tilde M=\Pi_{a=1}^r M_a\to \Q_\e^N\subset \Ee^{N+1},$$
 with $f_a$ spherical or $M_a$ one-dimensional for every $1\leq a\leq r$, a parallel vector bundle isometry 
$\phi\colon \tilde M\times \Ee^s \to N_{\tilde f}\tilde M$ and a local isometry 
$f_0\colon M_0\to (e_1+\Omega(\tilde f; \phi))\cap \Q_\e^{s-1}\subset \Ee^s$,   where $\phi_{\tilde p}(e_1)=\tilde f(\tilde p)$ for any $\tilde p\in \tilde M$, such that $f$ is given by (\ref{eq:def2}).
\end{corollary}

Arguing as in the proof of Theorem \ref{prop:n1b}, the general version of N\"olker's theorem follows easily from Theorems \ref{prop:n1c} and \ref{prop:n1e}.

\begin{theorem}\po\label{prop:n1d} Any isometric immersion  $f\colon M^n=M_0\times_{\rho} \Pi_{a=1}^r M_a\to \Q_\e^N$ of a warped product manifold  
whose second fundamental form is adapted to the product net is a warped product of immersions.
\end{theorem}

 \subsection[An application]{An application} 

  Let $f\colon\,M^n\to \Q_\e^N$ be an isometric immersion  with flat normal bundle. Then, it is well-known that for each point $p\in M^n$  there exist an integer $s(p)$  and unique principal  normals $\eta_1,\ldots,\eta_s\in N_fM(p)$  such that the tangent space splits orthogonally as
$$
T_pM=E_{\eta_1}(p)\oplus\cdots\oplus E_{\eta_s}(p)
$$
and the second fundamental form of $f$ splits accordingly as 
\be\label{eq:sff3}\alpha(X,Y)=\sum_{a=1}^s \<X^a,Y^a\>\eta_a,\ee
where $X^a$ denotes the $E_{\eta_a}$--component of $X$ for $1\leq a\leq s$. Assume, in addition,  that $s=s(p)$ is constant on $M^n$. In this case,  it is also well-known that  
each $\eta_a$ is smooth,
the dimension of $E_a:=E_{\eta_a}$ is constant, and $\E^f=(E_{1},\ldots, E_{s})$ is an orthogonal net  on $M^n$ such that  $E_{a}$ is umbilical for $1\leq a\leq s$,  and in fact spherical if $\rank\, E_{a}\geq 2$. Furthermore, if  $E_{a}$ is spherical then the restriction of $f$ to each leaf of $E_{a}$ is an spherical isometric immersion into $\Q_\e^N$ (see, e.g., \cite{dn}).

We have the following generalization of  Corollary \ref{cor3:ptubes} and its counterpart for the sphere and hyperbolic space.


\begin{corollary}\po \label{prop:fnb} Under the above assumptions, suppose there exists $1\leq r\leq s-1$ such that $E^\perp_{a}$ is totally geodesic for $1\leq a\leq r$.  Set $E_0:=\cap_{a=1}^r E^\perp_{a}.$
Then there exists locally a product representation $\Phi\colon\,  \Pi_{i=0}^{r} M_i\to M$ of $(E_0, E_1,\ldots, E_r)$, which is an isometry with respect to a quasi-warped product metric on $\Pi_{i=0}^{r} M_i$, such that $f\circ \Phi$ is given as described before Theorems  \ref{prop:n1c} or   \ref{prop:n1e}, according as $\epsilon=0$ or $\epsilon\in \{-1, 1\}$, respectively.
\end{corollary}

\proof That there exists locally a product representation $\Phi\colon\,  \Pi_{i=0}^{r} M_i\to M$ of the net $(E_0, E_1, \ldots, E_r)$ follows from Proposition \ref{prop:dec}. Since $E_{a}$ is umbilical and $E^\perp_{a}$ is totally geodesic for $1\leq a\leq r$, Proposition \ref{cor:prodwarp} implies that the metric induced by $\Phi$ is a quasi-warped product metric. Finally, the second fundamental form of $f$ is adapted to $(E_0, E_1, \ldots, E_r)$ by (\ref{eq:sff3}), hence the second fundamental form of $f\circ \Phi$ is also  adapted to the product net of  
$\Pi_{i=0}^{r} M_i$. Thus either Theorem \ref{prop:n1c} or  Theorem  \ref{prop:n1e}  applies to $f\circ \Phi$,
according as $\epsilon=0$ or $\epsilon\in \{-1, 1\}$, respectively.\vspace{2ex}\qed

Finally, given  $1 \leq \ell\leq s$,   below we give  some known conditions for $E^\perp_{\eta_\ell}$ to be totally geodesic (cf. \cite{dn}):

\begin{proposition}\po The following holds:
\begin{itemize}
\item[$(i)$] $E^\perp_{\eta_\ell}$ is totally geodesic if and only if it is integrable and $\eta_a$ is parallel along $E_{\eta_\ell}$ for every $1\leq a\leq s$ with $a\neq \ell$.
\item[$(ii)$] If the vectors $\eta_a-\eta_{\ell}$ and $\eta_b-\eta_{\ell}$
are everywhere linearly independent for any pair of indices
\mbox{$1\leq a\neq b\leq s$} with $a,b\neq \ell$ then $E^\perp_{\eta_\ell}$ is integrable.
\end{itemize}
\end{proposition}
\proof The Codazzi equations  yield
\be\label{eq:codz}
\nap_{X_\ell}\eta_a=\<\nabla_{X_a}X_a,X_\ell\>(\eta_a-\eta_\ell)
\;\;\;\mbox{if}\;\;\;a\neq \ell,
\ee
and
\be \label{eq:codz2}
\<\nabla_{X_b}X_a,X_\ell\> (\eta_a-\eta_{\ell})
=\<\nabla_{X_a}X_b,X_\ell\> (\eta_b-\eta_{\ell})\;\;\;\mbox{if}\;\;\;a, b\neq \ell
\ee
for all unit vector fields
$X_b\in \Gamma(E_{\eta_b})$, $X_a\in \Gamma(E_{\eta_a})$ and
$X_{\ell}\in \Gamma(E_{\eta_{\ell}})$.

Equation (\ref{eq:codz}) implies that $\eta_a$ is parallel along $E_{\eta_\ell}$ for every $1\leq a\leq s$ with $a\neq \ell$ if 
$E^\perp_{\eta_\ell}$ is totally geodesic. Conversely, if $E^\perp_{\eta_\ell}$ is integrable then the first factors in both sides 
of (\ref{eq:codz2}) coincide, hence both must vanish if $a\neq b$. On the other hand, from (\ref{eq:codz})  and the assumption that $\eta_a$ is parallel along $E_{\eta_\ell}$ for every $1\leq a\leq s$ with $a\neq \ell$ we obtain that $\<\nabla_{X_a}X_a,X_\ell\>=0$ for every $1\leq a\leq s$ with $a\neq \ell$. Thus $E^\perp_{\eta_\ell}$ is totally geodesic. The assertion in $(ii)$ follows immediately from 
(\ref{eq:codz2}).\qed

{\renewcommand{\baselinestretch}{1}
\hspace*{-25ex}\begin{tabbing}
\indent  \=  Universidade Federal de S\~ao Carlos \\
\indent  \= Via Washington Luiz km 235 \\
\> 13565-905 -- S\~ao Carlos -- Brazil \\
\> e-mail: tojeiro@dm.ufscar.br \\
\end{tabbing}}

\end{document}